\documentclass[a4paper,11pt]{article}

\usepackage{amsmath,amsfonts,amssymb}
\usepackage{ stmaryrd }
\usepackage{graphics}
\usepackage{enumerate}
\usepackage{mathrsfs}
\usepackage[amsmath,thref, thmmarks]{ntheorem}
\usepackage{mathpartir}
\newtheorem*{theorem-non}{Theorem}

\input{diagxy}
\xyoption{curve}

\newbox\anglebox % large pullback angle
\setbox\anglebox=\hbox{\xy \POS(75,0)\ar@{-} (0,0) \ar@{-} (75,75)\endxy}
 
\newbox\angleboxr % reverse large pullback angle
\setbox\angleboxr=\hbox{\xy \POS(0,0)\ar@{-} (0,75) \ar@{-} (75,0)\endxy}
 
\newbox\sanglebox % small pullback angle
\setbox\sanglebox=\hbox{\xy \POS(50,0)\ar@{-} (0,0) \ar@{-} (50,50)\endxy}
 
\newbox\sangleboxr % small reverse pullback angle
\setbox\sangleboxr=\hbox{\xy \POS(0,0)\ar@{-} (0,50) \ar@{-} (50,0)\endxy}
 
\newbox\sangleboxf % small flipped pullback angle
\setbox\sangleboxf=\hbox{\xy \POS(50,50)\ar@{-} (50,0) \ar@{-} (0,50)\endxy}
 
\newbox\angleboxf % flipped pullback angle
\setbox\angleboxf=\hbox{\xy \POS(75,75)\ar@{-} (75,0) \ar@{-} (0,75)\endxy}

\newbox\sangleboxfr % small flipped reverse pullback angle
\setbox\sangleboxfr=\hbox{\xy \POS(0,50)\ar@{-} (50,50) \ar@{-} (0,0)\endxy}
 
\newbox\angleboxfr % small flipped reverse pullback angle
\setbox\angleboxfr=\hbox{\xy \POS(0,75)\ar@{-} (75,75) \ar@{-} (0,0)\endxy}

\newcommand{\moddiag}[1]{\ensuremath{\operatorname{MDiag}\!\left({#1}\right)}}
\newcommand{\strdiag}[1]{\ensuremath{\operatorname{Diag}\!\left({#1}\right)}}
\newcommand{\evalmod}[2]{\ensuremath{\operatorname{Ev}_{#1}\!\left({#2}\right)}}
\newcommand{\pair}[1]{\ensuremath{\langle {#1} \rangle}}
\newcommand{\sublat}[2]{\ensuremath{\operatorname{Sub}_{#1}\!\left({#2}\right)}}
\newcommand{\mo}[1]{\ensuremath{\operatorname{#1}}}
\newcommand{\mb}{\ensuremath{\mathbin}}

\newcommand{\cat}[1]{\ensuremath{\mathcal{#1}}}
\newcommand{\thry}[1]{\ensuremath{\mathbb{#1}}}
\newcommand{\synt}[2]{\ensuremath{\mathcal{#1}_{{#2}}}}
\newcommand{\alg}[1]{\ensuremath{\mathbf{#1}}}

\newcommand{\modcat}[1]{\ensuremath{\operatorname{Mod}(#1)}}
\newcommand{\strcat}[1]{\ensuremath{\mathrm{Str}({#1}})}

\newcommand{\topo}[1]{\ensuremath{\mathscr{#1}}}
\newcommand{\Sets}{\ensuremath{\mathbf{Set}}}

\newcommand{\op}[1]{\ensuremath{{#1}^\mathrm{op}}}

\newcommand{\Sh}[1]{\protect\ensuremath{\operatorname{Sh}\!\left(#1\right)}}

\newcommand{\sem}[1]{\ensuremath{[\![{#1}]\!]}}
\newcommand{\csem}[2]{\ensuremath{[\![{#1}\mb|{#2}]\!]}}% withcontext

\newcommand{\alle}[1]{\forall {#1}\mathpunct .}

\newcommand{\fins}[1]{\exists {#1}\mathpunct .}

\newcommand{\theory}{\ensuremath{\mathbb{T}}}
\newcommand{\cterm}[2]{\ensuremath{\left \{ {#1}\ \; \vrule \; \ {#2}\right \}}}
\newcommand{\syntob}[2]{\ensuremath{[{#1}\;|\;{#2}]}}
\newcommand{\funktor}[3]{\ensuremath{\operatorname{#1} \! :\!  #2 \to<125> #3}}
\newcommand{\FIC}[2]{\ensuremath{{#1}.{#2}}}
\newcommand{\Image}[1]{\ensuremath{\operatorname{Im}(#1)}}
\title{Constructive completeness and non-discrete languages}
\author{Henrik Forssell and Christian Esp\'indola}
\date{April 11, 2016}

\begin{document}
\maketitle
%\tableofcontents

%ABSTRACT
\begin{abstract}
We give an analysis and generalizations  of some long-established constructive completeness results in terms of categorical logic and pre-sheaf and sheaf semantics.
% and add some new results in this area which flow from the analysis. 
The purpose is in no small part conceptual and organizational: from a few basic ingredients arises a more unified picture connecting constructive completeness with respect to Tarski semantics, to the extent that it is available, with various completeness theorems in terms of presheaf and sheaf semantics (and thus with Kripke and Beth semantics). From this picture are obtained both (``reverse mathematical'') equivalence results and new constructive completeness theorems; in particular, the basic set-up is flexible enough to obtain strong constructive completeness results for languages of arbitrary size and languages for which equality between the elements of the signature is not decidable. 
%
%That is, without the  assumption that equality of non-logical symbols in the language can be decided.  

\end{abstract}
%

%INTRODUCTION
\section{Introduction}

Starting with the G\"odel-Kreisel theorem, it has long been well known that the classically ``standard'' semantics---Tarski structures for classical first-order logic (FOL), and Kripke or Beth structures for intuitionistic FOL---are insufficient in a constructive metatheory. For instance, the assumption of strong completeness for intuitionistic FOL with respect to Tarski,  Kripke, or Beth semantics in a metatheory such as \textbf{IZF}, \textbf{HAS}, or \textbf{HAA}\footnote{Intuitionistic Zermelo Fraenkel set theory; full intuitionistic second-order arithmetic; intuitionistic second-order arithmetic with arithmetic comprehension}  implies the law of excluded middle (LEM), while weak completeness implies Markov's principle (MP) (see \cite{mccarty:08}, \cite{kreisel:62}, \cite{mccarty:94}). On the other hand,  constructive completeness theorems exist e.g.\ with respect to formal space valued models,   and in sheaf toposes more generally (see e.g.\ \cite{MP}, \cite{coquandsmith:95}, \cite{elephant1}). In a sense intermediate between sheaf semantics and the standard semantics of Tarski and Kripke, completeness was also shown (albeit assuming the Fan Theorem) to hold for ``countable'' intuitionistic first-order theories with respect to fallible Kripke and Beth semantics by Veldman \cite{veldman:76} and by de Swart \cite{deswart:76}.% in which of nodes that force $\bot$ are allowed.

The insufficiency of the standard semantics can be taken to suggest that constructive model theory should instead be carried out with respect to sheaf models, or some variant thereof. On the other hand, the theorems of Veldman and de Swart indicate that ordinary Tarski semantics can have an important role to play. Thus the underlying conceptual and motivating question of this paper is 
%
%what the ``correct'' constructive notion of a model is, and 
%
the role of Tarski semantics in constructive model theory, or how much ``mileage'' one can get out of ordinary Tarski-models in a constructive setting. 

Developments of logic in a constructive setting usually assume that the signatures are, if not in some sense countable, then at least discrete; i.e.\ that decidability (LEM) holds for equality   of the basic function and relation symbols. This precludes various classical constructions, such as adding the elements of an arbitrary domain as constants to the language. Or similarly in categorical logic, constructing the internal language of an arbitrary category. Here, we drop this restriction.  

In a classical metatheory (with choice), the link between  Tarski completeness for classical first-order logic (or its so-called coherent fragment) and Kripke completeness for intuitionistic first-order logic is put to light in a theorem  attributed to A.\ Joyal in \cite{makkaireyes}.
The theorem involves, among other things, the technique of considering a theory in terms of an ``approximation'' in a weaker fragment, known sometimes, or in some cases, as \emph{Morleyization}. 
However, the theorem is not immediately applicable   in a constructive setting as it relies on the assumption of Tarski completeness for the coherent fragment of FOL. The proof (\cite[Thm. 6.3.5]{makkaireyes})  also uses classical techniques.
%The latter  (and Kripke completeness for FOL) fails, of course, in general, to hold in a constructive meta-theory. 
Nevertheless, while Tarski completeness for the coherent fragment fails to hold in a constructive meta-theory, constructive completeness results for less expressive fragments of FOL exist.  By giving a constructive formulation and proof of Joyal's theorem, and using the same ``approximation'' technique, these can be exploited to give completeness results for stronger fragments of FOL, and for FOL itself, in suitable pre-sheaf and sheaf toposes. For instance, the aforementioned completeness theorems for fallible Kripke and Beth semantics can be recovered (and given  new proofs) in this way. The purpose and aim of this paper is in no small part to ``tell this story''; that is, to give a unified and conceptual account connecting constructive Tarski completeness, to the extent that it exists, with completeness results in terms of traditionally studied Kripke and Beth-style models. And, furthermore,  to give this account using essentially only the two basic ingredients of Joyal's theorem and of Morleyization to fragments for which Tarski completeness holds. In addition to displaying the connections between already established results, such a clearer conceptual picture can also serve to suggest further, and new, ones. As an instance of this, we extend constructive completeness results to languages and theories that are not enumerable, and for which  equality of non-logical symbols in the language can not be decided. (We also do not, in general, assume that  the sentence $\fins{x}x=x$ is valid). In particular, we give a strong completeness theorem for the disjunctive-free fragment of FOL, over such languages, with respect to fallible Kripke models, and a strong completeness theorem for full FOL with respect to fallible, ``generalized'' Beth models. We also draw attention to a conceptual link between Beth semantics and a ``least coverage forcing the correct interpretation of disjunctions'',   and  give a constructive version and proof of the completeness theorem of \cite{gabbay:77} with respect to Beth models in which only the forcing clause for disjunction is used.
The paper is structured into the following parts. Section \ref{section: mod tarsk completeness} contains preliminaries and notes on notation and terminology, and recapitulates Tarski completeness for the regular fragment and enumerable coherent theories. The regular completeness theorem for arbitrary theories and signatures is more fully presented in  \cite{forsselllumsdaine}. The completeness theorem for enumerable positive coherent theories is known, but is included  for conceptual self-containement.   (We also note the equivalence between this theorem and the Fan theorem, 
and draw a corollary concerning completeness for classical first-order theories.)
% (\ref{corollary: completeness of classical fo theories}).
%
Section \ref{section: joyals theorem} gives a constructive reformulation and proof of Joyal's theorem. 
%
%(\ref{theorem: joyals theorem})
%
Section \ref{section: sheaf completeness} introduces and analyses certain coverages on the category of models, and gives a covering lemma which allows Joyal's theorem to be stated with respect to a poset of structures and  homomorphic inclusions. This then linked with the instances of Tarski completeness in Section \ref{section: mod tarsk completeness} to give the aforementioned Kripke and Beth completeness results
%, both for enumerable languages and theories, and for arbitrary theories over non-discrete languages. 
theories over non-discrete languages. The specialization to the constructive version of the theorem of \cite{gabbay:77} marks the end.

\section{Constructive Tarski completeness}\label{section: mod tarsk completeness}%\section{Joyal's theorem}\label{section: Joyal's theorem}

%We state Joyal's theorem in a form which is convenient for our purposes and give a constructive presentation of the proof. We %also introduce in this section notation and definitions that will be used further down.

\subsection{Preliminiaries}\label{subsection: preliminaries}

\subsubsection{Theories, models and diagrams}
\label{subsection: models and diagrams}
\label{subsubsection: models and diagrams}

In what follows we shall fix our metatheory to be \textbf{IZF} and simply use ``constructive'' to mean that we are working in this setting. Correspondingly, ``classically'' means in the metatheory \alg{ZFC}.
We fix the following terminology, conventions, and notation. 
By ``finite'' and ``countable'' we mean cardinal finite and  isomorphic to \thry{N}, respectively. A ``list'' is a finite list. \Image{\alg{a}} denotes the set of elements of the list \alg{a}, and  $l(\alg{a})$ its length. When deemed safe, we shorten $\alg{a}\in A^{{l(\alg{a})}}$ to $\alg{a}\in A$.
 A subset of a set $A$ is decidable if it is given in terms of a function $f\! :\!A\to<100>2$ as  \cterm{a\in A}{f(a)=1} and semi-decidable if it is given in terms of a function   $f\! :\!A\to<100>2^{\thry{N}}$ as  \cterm{a\in A}{\fins{n}(f(a))(n)=1}. By an \emph{enumerable}\ set we mean a semi-decidable subset of a (perhaps implicit) countable set.  

Let $\Sigma$ be a first-order signature. We generally assume that $\Sigma$ is \emph{relational} in the sense that it has no function symbols. Thus all functions and constants are taken to be represented by relations and appropriate axioms over $\Sigma$.
%We nevertheless use constant symbols when considering theories of models and related constructions. But this should then be considered a presentantional shortcut, relying on the existence of a provability-preserving and reflecting translation into a theory with relation symbols only, where the constants are represented by unary predicates (cf.\ \cite{bell-machover:1977}).  
%
% It is mostly function symbols with positive arity we wish to avoid, however, and for  convenience we do use constants e.g.\ when dealing with diagrams and theories of models, instead of explicitly rewriting them as unary relations and adding appropriate axioms.
 Furthermore, we  assume,  that $\Sigma$ is single-sorted. 
%
%It is straightforward to rewrite what follows to the many-sorted case.
%
%Thus $\Sigma$ is simply a set of relation symbols with associated (finite) arities. We do not assume, unless otherwise stated, that $\Sigma$ is {discrete}---i.e.\ that we can decide whether the symbols in $\Sigma$ are the same or not---nor that it is {decidable}---i.e.\ that we can decide whether a  symbol is in the signature or not. (As usual, however, the logical symbols, including variables, are both discrete and decidable, in this sense). 
%
Thus $\Sigma$ is simply a set of relation symbols with associated (finite) arities. We do not assume, unless otherwise stated, that $\Sigma$ is \emph{discrete}---i.e.\ that LEM holds for equality betweeen the elements of  $\Sigma$. (As usual, however, the logical symbols, including variables, are discrete, and disjoint from $\Sigma$). 
Following \cite[D1]{elephant1}, we consider theories in FOL and fragments of FOL formulated in terms of sequents of the form  $\phi\vdash_{\mathbf{x}}\psi$.  These can be read as $\alle{\alg{x}}\phi\rightarrow \psi$. 
The list of variables \alg{x} is required to be a \emph{context} for both $\phi$ and $\psi$ in the sense that it is a list of distinct variables containing (at least) the free variables of the formula (see \cite[D1.1.4]{elephant1}). We write \FIC{\alg{x}}{\phi} for a \emph{formula-in-context}. We write \syntob{\alg{x}}{\phi} for a formula-in-context identified up to $\alpha$-equivalence. (That is, up to renaming of bound variables and variables in the context, as in \cite[D1.4]{elephant1}.) The context, in both cases, is \emph{canonical} if it contains only the free variables of the formula, listed in order of first-appearance. 
The logic is ``free'', in the sense that the sequent $\top\vdash \fins{x}x=x$ is, in general, not derivable. 
 The main fragments we shall be referring to are: the \emph{Horn} fragment\footnote{To prevent confusion with otherwise standard usage of ``Horn clause'' and ``Horn formula'', note the usage here (following \cite{elephant1}) of ``Horn formula'' as simply a formula which is a conjunction of atomic formulas, and ``Horn sequent'' as simply a sequent with such Horn formulas as antecedent and consequent. }, consisting of sequents with formulas over $\Sigma$ involving only the logical constants $\top$ and $\wedge$; the \emph{regular} fragment $\top, \wedge, \exists$; the $\mathit{regular_{\bot}}$ fragment $\top, \wedge, \exists$, and $\bot$; the \emph{positive coherent} fragment $\top, \wedge, \vee$, and $\exists$; the $\mathit{coherent}$ fragment $\top, \wedge, \vee, \exists$, and $\bot$; and, of course, full FOL. 
 %the $\vee$-$\mathit{free}$ fragment involving $\top, \wedge, \bot, \rightarrow, \exists$, and $\forall$. 
 Deduction rules and further details can be found in \cite[D1.3]{elephant1}. The distinguished relation symbol $=$ of equality is included in all languages under consideration. If a theory \theory\ proves a sequent $\phi\vdash_{\mathbf{x}}\psi$ we sometimes write $\phi\vdash_{\mathbf{x}}^{\theory}\psi$ instead of  $\theory \vdash (\phi\vdash_{\mathbf{x}}\psi)$. Provable in the empty theory is then written $\phi\vdash_{\mathbf{x}}^{\emptyset}\psi$.

One would usually say that a theory is regular, for instance, if it is axiomatizable by regular sequents. Thus a Horn theory would also be a regular theory etc. For brevity, however, we also mean to indicate what fragment we are considering when we say that a theory is this or that. Thus when we say e.g.\ that \theory\ is a coherent theory and  $\phi$ is a formula of \theory, we mean, in particular,  that $\phi$ a formula in the coherent fragment over the signature of \theory. In the same vein, if we say that a theory is discrete, we mean that it is over a discrete signature. If we say that a theory is enumerable\ we mean both that it is over a enumerable\ signature and that the set of axioms is enumerable

We say that a coherent sequent is on \emph{canonical form} if it is on the form $\phi\vdash_{\mathbf{x}}\exists{\mathbf{y}_0}\psi_0\vee\ldots\vee\exists{\mathbf{y}_n}\psi_n$ where $\phi$ and all $\psi_i$ are Horn formulas, or on the form $\phi\vdash_{\mathbf{x}}\bot$ where $\phi$ is Horn.  A regular sequent is on canonical form if it is so as a coherent sequent. Every coherent ($\mathrm{regular}_{\bot}$, regular) theory can be axiomatized by coherent ($\mathrm{regular}_{\bot}$, regular)   sequents on canonical form (see \emph{ibid.}), and  we assume that they are.
%, and that the truth of general sequents can be defined in terms of truth of sequents on canonical form.
%
%We say that a signature $\Sigma$ (or theory over $\Sigma$) is discrete if the set of non-logical symbols is discrete, and we say that  a theory is enumerable\ if both the set of relation symbols and the set of axioms are enumerable\ sets. 

By \emph{ Tarski structure} for a signature, and \emph{Tarski model} for a theory, we mean the usual notion of a domain set with interpretations of the relation symbols in terms of subsets, and the interpretations of the connectives $\bot$, $\top$, $\wedge$, $\vee$, and $\exists$ by the usual set-theoretic interpretations (as well as $\rightarrow$ and $\forall$, but we do not actually consider Tarski models for anything above the coherent fragment, with the exception of Corollary \ref{corollary: completeness of classical fo theories}). 
The domain need not be inhabited or non-empty. 
As for the interpretation of the  equality relation, we reserve ``structure'' and ``model'' for the case where equality is interpreted as the identity relation, and use \emph{diagram} and \emph{model diagram} for the case where equality is interpreted as a congruence relation\footnote{We avoid using ``diagram'' when extending a language with a structure, talking instead of the language and theory of the structure.}. That is to say, a diagram for a relational signature $\Sigma$ consists of a $\Sigma$-structure \alg{M}  together with an equivalence relation $E$ on $|\alg{M}|$ which respects the relations interpreting the symbols of $\Sigma$. The interpretation  $\csem{\alg{x}}{\phi}^{\alg{M}}$ of a formula-in-context \FIC{\alg{x}}{\phi} is then defined in the usual way, but interpreting $=$ as $E$.

In part to notationally distinguish diagrams, in this sense, from structures, we consider and write a diagram \alg{M} as a pair $(D,F)$ where $D=|\alg{M}|$ is the domain and $F=\cterm{\pair{\syntob{\alg{x}}{\phi},\alg{d}}}{\phi\ \textnormal{is Horn, and}\ \alg{d}\in\csem{\alg{x}}{\phi}^{\alg{M}} }$, where $\csem{\alg{x}}{\phi}^{\alg{M}}$ is the extension of \FIC{\alg{x}}{\phi} in \alg{M}. 
We refer to an element of $F$ as a \emph{fact}. 
 For a signature $\Sigma$ and theory \theory\ we write \strcat{\Sigma} and   \modcat{\theory} for the category of structures and homomorphisms and models and homomorphisms, respectively.  We write \strdiag{\Sigma} and \moddiag{\theory} for the category of diagrams and homomorphisms and  model diagrams and homomorphisms, respectively, where a homomorphism $h:(D_1, F_1)\to<150>(D_2, F_2)$ between diagrams is a left-total relation $h\subseteq D_1\times D_2$ such that 
\begin{enumerate}
\item $h(d_1,d_2)\wedge\pair{\syntob{x,y}{x=y},d_2,d'_2}\in F_2\rightarrow h(d_1,d'_2)$; and
\item $\pair{\syntob{\alg{x}}{\phi},\alg{d}_1}\in F_1)\wedge h(\alg{d}_1,\alg{d}_2)\rightarrow (\pair{\syntob{\alg{x}}{\phi},\alg{d}_2}\in F_2)$, for all (atomic) $\mathrm{Horn}$ formulas $\phi$ over $\Sigma$ and all  $\alg{d}_1\in D_1$ and $\alg{d}_2\in D_2$. 
\end{enumerate}
(Here $h(\alg{d}_1,\alg{d}_2)$ stands for the expected conjunction.  As further notational shortcuts, we allow ourselves to use function notation for homomorphisms between diagrams when no confusion threatens. For instance we might write $\phi[h(d_1)/x]$ instead of $\alle{d_2}h(d_1,d_2) \rightarrow \phi[d_2/x]$. We sometimes write $\phi[\alg{d}/\alg{x}]\in F$, or simply $\phi[\alg{d}]\in F$, instead of $\pair{\alg{d},\syntob{\alg{x}}{\phi}}\in F$, for brevity.) 
%
%For a relational signature $\Sigma$ 
Then there is an adjoint equivalence 
\[\bfig \morphism/{@{>}@/^1em/}/<1000,0>[\strcat{\Sigma}`\strdiag{\Sigma};i] 
\morphism|b|/{@{<-}@/_1em/}/<1000,0>[\strcat{\Sigma}`\strdiag{\Sigma};q]
\place(500,0)[\simeq]  \efig\]
where $i$ is the inclusion and $q$ is by taking quotients, in the expected way. The unit \funktor{k}{(D,F)}{q(D,F)} preserves and reflects the intepretation of formulas and the truth of sequents. Thus, in particular, the equivalence restricts to models,   $\modcat{\theory}\simeq \moddiag{\theory}$.

The diagram notation may seem cumbersome at first, but it is convenient for working with presentations of structures. Let a \emph{presentation} of a diagram, or \emph{pre-diagram}, be a pair $(D,F)$ where $D$ is a set and $F$ is set of  facts over $\Sigma$ and $D$, in the sense above; that is to say, $F$ is  a set of pairs \pair{\syntob{\alg{x}}{\phi},\alg{d}} with $\phi$ Horn and $\alg{d}\in D^{l(\alg{x})}$. For a pre-diagram $(D,F)$ the least diagram containing it is the diagram \emph{generated by} $(D,F)$. A homomorphism of pre-diagrams is a left-total relation that is a homomorphism of the generated diagrams.  The satisfaction relation $\vDash$ for pre-diagrams is defined as satisfaction in the generated diagram.

For  $(D,F)$ is a pre-diagram, the theory $\thry{D}_{(D,F)}$ of $(D,F)$ is defined as expected by extending $\Sigma$ with  $D$ as constants (or unary predicates, see below) and letting \cterm{\top\vdash \phi[\alg{d}/\alg{x}]}{\pair{\syntob{\alg{x}}{\phi}, \alg{d}}\in F} be axioms. The diagram generated by $(D,F)$ can then be defined as $(D,\cterm{\pair{\syntob{\alg{x}}{\phi}, \alg{d}}}{\top\vdash^{\thry{D}_{(D,F)}} \phi[\alg{d}/\alg{x}]})$. If \theory\ is a theory over $\Sigma$ we write $\theory_{(D,F)}$ for the union of \theory\  and $\thry{D}_{(D,F)}$. 

When constants are not discrete, replacing constants with variables in proofs becomes more problematic; that is, one cannot in general replace the same constant with the same variable throughout a formula. This makes the interplay between a theory and the theory of one of its diagrams a little more intricate.   The following lemmas, which are completely straightforward for discrete signatures, are shown also to hold for non-discrete signatures in  \cite{forsselllumsdaine} and stated here for reference.
\begin{lemma}\label{lemma: replacing constants with variables}
Let \theory\ be a
%
%\textup{(}first-order\textup{)}
%
theory over $\Sigma$. Let $C$ be a set of constants disjoint from $\Sigma$, and write $\Sigma^C=\Sigma\cup C$. Suppose $\phi\vdash_{\alg{x}}\psi$ is a  first-order sequent over $\Sigma^C$ which is provable from axioms in \theory. Then there exists a sequent $\phi'\vdash_{\alg{x},\alg{y}}\psi'$ over $\Sigma$ and a ``valuation'' function $\funktor{f}{\alg{y}}{C}$ such that;   
\begin{enumerate}[(i)]
\item $\phi'\vdash_{\alg{x},\alg{y}}\psi'$ is provable from \textup{(}the same\textup{)} axioms in \theory; and
\item $\phi'[f]=\phi$, and $\psi'[f]=\psi$. 
\end{enumerate}
\end{lemma}
\begin{lemma}\label{lemma: getting rid of constants from the model}
Let \theory\ be a 
%
%\textup{(}first-order\textup{)}
%
 theory over $\Sigma$ and $(D,F)$ be a \textup{(}pre-\textup{)}diagram. Let $\FIC{\alg{x}}{\psi}$ and $\FIC{\alg{x},\alg{y}}{\phi}$ be  first-order $\Sigma$-formulas-in-context. Let  \alg{c} be a tuple of elements in $D$. Suppose $\theory_{D,F}$ proves the sequent $\phi[\alg{c}/\alg{y}]\vdash_{\alg{x}}\psi$. Then there is a regular formula $\xi$ in context \alg{y} in $\Sigma$ such that $(D,F)\vDash\xi[\alg{c}/\alg{x}]$ and \theory\ proves the sequent $\xi\wedge\phi\vdash_{\alg{x},\alg{y}}\psi$.
\end{lemma}
%
%Note that for coherent \FIC{\alg{x}}{\phi} and $\alg{d}\in D$, $(D,F)\vDash \phi[\alg{d}/\alg{x}]$ if and only if $\top\vdash^{\thry{D}_{D,F}} \phi[\alg{d}/\alg{x}]$.

As mentioned in the beginning, we generally assume that signatures are relational  and thus that any language with function symbols has been translated into an equivalent one without them (cf.\ \cite{bell-machover:1977}). This is in principle so also when extending a signature with the language of one of its diagrams, as above. In practice this becomes burdensome, however, and we leave the translation  implicit and swept under the rug.  

 Finally, we say that a diagram $(D,F)$ is a \emph{subdiagram} of $(D',F')$ if $D\subseteq D'$ and $F\subseteq F'$. We write $(D,F)\subseteq (D',F')$.  Note that the inclusion $D\subseteq D'$ induces a homomorphism   $\funktor{i}{(D,F)}{(D',F')}$ by $i(d,d')\Leftrightarrow (d=d')\in F'$.  (This homomorphism need not be a monomorphism.)
%
%We say that a diagram $(D,F)$ is finite, discrete, or enumerable\ when both $D$ and $F$ are finite, discrete, or enumerable, respectively.
%
A  diagram is \emph{finite} if the domain is finite and the interpretations of $=$ and all relation symbols are finite. For enumerable\ signature $\Sigma$, a diagram is \emph{enumerable}\ if the domain is enumerable\ and  the interpretations of $=$ and all relation symbols are enumerable.  For discrete signature $\Sigma$, a diagram is \emph{discrete} if the domain is discrete. 
We define a \emph{bounded diagram} to be a triple $(D,F,n)$ where $(D,F)$ is a diagram, $n\in \thry{N}$, and the elements of the domain $D$  are pairs where the second component is a natural number less than or equal to $n$. 
Clearly, any diagram is canonically isomorphic to a bounded one (with bound $0$, say).
%
%We use this as a device   to be able to ``freely'' add elements to diagrams and recognize an element as added. 
We shall in fact mostly restrict to bounded diagrams, but leave the bound $n$ notationally implicit. Let $\operatorname{Diag}_b(\Sigma)$ be the category of bounded diagrams and diagram homomorphisms.

\subsubsection{Syntactic categories, Morleyization, and exploding models}\label{subsubsection: morleyization} %%%%%%%%%%%%%%%%%%%%%%%%%%%%%%
%%%%%%%%%%%%%%%%%%%%%%%%%%%%%%%%%%%%%%%%%%%%%%%%%%%%%%%%%%%%%%%%%%%%%%%%%%%%%%%%%%%%%%%%%%%%%%%%%%%%%%%%%%%%%%%%%
%

Recall from e.g.\ \cite[D1.4]{elephant1} the \emph{syntactic category} \synt{C}{\theory} of a theory \theory, consisting of formulas-in-context of the language of the theory. For coherent \theory, the category   \synt{C}{\theory} is a coherent category and models of \theory\ in a coherent category \cat{D} can be considered as coherent functors $\synt{C}{\theory}\to<125>\cat{D}$. Similarly for e.g.\ regular and first-order theories, see \emph{loc.\ cit.} for precise statements and details. 
The functor $\synt{C}{\theory}\to<125>\cat{D}$ is \emph{conservative} (see e.g.\ \cite{elephant1}) if and only if the corresponding model is conservative, or complete; that is, if only provable sequents are true.

We refer to the rewriting of a theory into a theory of a less expressive fragment as ``Morleyizing'' the theory, after the rewriting of a classical first-order theory as an equivalent coherent theory (as in e.g.\ \cite[D1.5.13]{elephant1}). In categorical terms, the syntactic category \synt{C}{\theory} of, say, an intuitionistic first-order theory \theory\ is a Heyting category, and thus also a coherent and regular category. The Heyting category \synt{C}{\theory} therefore has an internal coherent theory $\theory^{coh}_{\cat{C}_{\theory}}$ and an internal regular one $\theory^{reg}_{\cat{C}_{\theory}}$. These theories are equivalent in the sense that their syntactic categories are equivalent.
\[\synt{C}{\theory}\simeq \synt{C}{\theory^{coh}_{\cat{C}_{\theory}}}\simeq \synt{C}{\theory^{reg}_{\cat{C}_{\theory}}}\]
The categories of models of these theories will in general not be the same (unless we are considering classical theories and their coherent Morleyzations, see loc.cit.). Nevertheless, considering the category of models of the Morleyized theory can be fruitful, not least when the more expressive theory has ``too few models'' (see \cite{makkai:95} for another example).

% We refer to these as ``the theory of coherent (regular) models of \theory'', or  its coherent (regular) ``Morleyization'' for short, and denote it $\theory^m$ (leaving to context rather than notation which fragment we are considering). The theory and its Morleyization have equivalent syntactic categories,
%
%\[\synt{C}{\theory}\simeq \synt{C}{\theory^m}\]
%
%If \theory\ is classical first-order, the coherent Morleyization will have the same models as \theory\ in all Boolean coherent categories (see \emph{loc.cit.}), and is in that sense equivalent to it. In general, when Morleyizing an intuitionistic theory to a coherent one, say, or a coherent to a regular, this is not the case. Nevertheless, considering the category of models of the Morleyized theory can be fruitful, not least when the more expressive theory has ``too few models'' (see \cite{makkai:95} for another example).   

There is some leeway concerning what precisely one takes the internal language and theory of a category to be. In our case we also start from a given theory and not from a given category. We therefore   write down explicitly what we shall take  the \emph{regular} and \emph{coherent Morleyizations} of a first-order theory  \theory\ over a signature $\Sigma$ to be. Other fragments are similar. 

%Other Morleyizations, such as the regular$_{\bot}$ Morleyization of a coherent theory, are  similar.
%
%In this section, we display the canonical context of a formula by writing $\phi(\alg{x})$ where \alg{x} is a possibly empty list of variables. We use square brackets for simultaneous substitution.
%
Let $\Sigma^m$ be the signature extending $\Sigma$ with, for each first-order formula ${\phi}$ over $\Sigma$, in canonical context \alg{x}, say, a relation symbol \textbf{P}$_{\phi}$ with arity the length of  $\alg{x}$. We write $P_{\phi}$ for the atomic formula over $\Sigma^m$ obtained by assigning $\alg{x}$ to the arity of \textbf{P}$_{\phi}$. 
Consider the following axioms.
\begin{itemize}
\item[(\emph{Thry})] \hfill \\ For every sequent $\phi\vdash_{\alg{x}}\psi$ provable in \theory, the axiom
\[P_{\phi}\vdash_{\alg{x}}P_{\psi}\]
%
%\item[(\emph{Subst})] \hfill \\ For formulas \FIC{\alg{x}}{\phi}, \FIC{\alg{y}}{\psi}, and list of terms \FIC{\alg{y}}{\alg{t}}, of the same length as \alg{x},  over $\Sigma$ such that $\psi=\phi[\alg{t}/\alg{x}]$ over $\Sigma$
%For every formula-in-context \FIC{\alg{x}}{\phi(\alg{x}')}, list of terms \alg{t} of the same lenght as \alg{x} and in context \alg{y}, and with $\phi[\alg{t}/\alg{x}']=\psi(\alg{y}')$, 
%
%the axiom
%  
%\[P_{\phi}[\alg{t}/\alg{x}]\dashv\vdash_{\alg{y}}P_{\psi}\] 
%
%where \alg{y} is the canonical context of $\psi$.
%
\item[(\emph{Atom})] \hfill \\ For every atomic formula $\phi$ over $\Sigma$ in canonical context \alg{x},   \[P_{\phi}\dashv\vdash_{\alg{x}}\phi\] 
\item[(\emph{True})] \hfill \\ \[P_{\top}\dashv\vdash_{}\top\]
\item[(\emph{Conj})] \hfill \\ For every conjunction $\theta=\phi\wedge\psi$ over $\Sigma$ in canonical context \alg{x}
 \[P_{\theta}\dashv\vdash_{\alg{x}}P_{\phi}\wedge P_{\psi}\]
\item[(\emph{Exist})] \hfill \\ For every existentially quantified formula $\theta=\fins{y}\phi$ over $\Sigma$ in canonical context \alg{x}
\[P_{\theta}\dashv\vdash_{\alg{x}}\fins{y}P_{\phi}\]
\item[(\emph{Disj})] \hfill \\ For every disjunction $\theta=\phi\vee\psi$ over $\Sigma$ in canonical context \alg{x}
 \[P_{\theta}\dashv\vdash_{\alg{x}}P_{\phi}\vee P_{\psi}\]
\item[(\emph{False})] \hfill \\ \[P_{\bot}\dashv\vdash_{}\bot\]
\end{itemize}
These axioms define the coherent Morleyization of \theory. The regular Morleyization is obtained by omitting the Disjunction axiom schema  and the False axiom. Notice that in, say, the regular Morleyization of a first-order theory, every regular formula is provably equivalent to an atomic formula; and that the sequent $\phi\vdash_{\alg{x}}\psi$ is provable in \theory\ if and only if  $P_{\phi}\vdash_{\alg{x}}P_{\psi}$ is provable in $\theory^m$. Notice further that if we Morleyize, say, a regular$_{\bot}$ theory \theory\ to a regular theory, then $\theory^m$ will prove $P_{\bot}\vdash_{\alg{x}}\phi$ for all regular \FIC{\alg{x}}{\phi} over $\Sigma^m$. A model diagram $(D,F)$ such that $P_{\bot}\in F$ will thus have all possible facts in $F$. The corresponding quotient will consist of a single point of which everything is true (it must be inhabited since $\theory^m$  proves $P_{\bot}\vdash \fins{x}\top$). Thus, or in that sense, it is an \emph{exploding} model diagram of \theory.

Finally we note  that the notion of enumerability for theories is not affected by Morleyization. We display this for reference.
%
%One\marginpar{NEW: added constructive content} can see Corollary \ref{corollary: classical joyals} as the %``constructive content'' of Joyal's Theorem as it appears in (\cite{makkaireyes}), making explicit that the non-%constructive part of that theorem is in the existence of such a category \cat{M}.
%
%
%\begin{definition}\label{Definition: s.d.}
%We say that a theory \theory\ over a signature $\Sigma$ is \emph{enumerable}\ (semi-decidable) if $\Sigma$ is a semi-%decidable subset of a countable set and (the set of axioms of) \theory\ is a semi-decidable subset of the set of %sequents over $\Sigma$.
%\end{definition}
%
\begin{lemma}
If \theory\ is enumerable\ then so is $\theory^m$.
\end{lemma}

%We can then conclude with the following:
%
%\begin{proposition}\label{proposition: Joyals corollary}
%Let \theory\ be a enumerable\ first-order ($\vee$-free) theory and let $\theory^m$ be its coherent (regular$_{\bot}$) Morleyization. Let  %$\modcat{\theory^m}$ and \topo{E}.
%\end{proposition}

 \subsection{Tarski completeness}\label{section: tarski completeness}%%%%%%%%%%%%%%%%%%%%%%%%%%%%%%%%%%%%%%%%%%%%%%%%%%%%%%%%%%

The dynamical method of \cite{costelombardiroy:01}, the chase algorithm of \cite{abiteboulhullvianu:95}, and similar methods\footnote{C.f.\ also \cite{adamekandrosicky:94} and \cite{makkai:90})} can be seen as simultaneous proof searches and (at least partial) completions of   structures to  models. In essence, one proceeds by repeatedly applying the axioms of the theory to the structure and adding the result; thus if, for instance, $\phi[a/x]$ is true in the structure and $\phi\vdash_x \fins{y}\psi$ is an axiom, one extends the domain with a fresh element $b$ and the interpretation in the least way  such that $\psi[a/x,b/y]$ is true. It is in several cases known, or at least folklore, that such methods can be used constructively to obtain completeness results for fragments of FOL. Although the object theories tend to be assumed countable or at least discrete. We summarize in Section \ref{Subsection: Regularbot theories} the relevant results from \cite{forsselllumsdaine} concerning the construction of a functorial  ``simultaneous chase'' to the case of regular theories with no size or discreteness constraints.  Section \ref{Subsection: Countable coherent theories and fan} displays 
the equivalence between the Fan theorem and  completeness for enumerable\ positive coherent theories with respect to enumerable\ model diagrams---equivalently of enumerable\  coherent theories with respect to ``possibly exploding'' or ``fallible'' enumerable\ model diagrams.
This equivalence can to a large extent be derived from the literature.  %notably from the work of Veldman's \cite{veldman:76} and Loeb's \cite{loeb:05}. 
In particular, Veldman's proof \cite{veldman:76} (which relies on the Fan Theorem)  of fallible Kripke completeness for first-order (enumerable) theories  implies also the  Tarski-completeness of positive coherent (enumerable) theories. Nevertheless, since we are, conceptually,  regarding first-order fallible Kripke completeness as flowing from the Tarski completeness of positive coherent theories, we supply a direct proof of the latter using the Fan Theorem. The converse, that this completeness theorem implies the Fan Theorem, is rather immediate, and we include a very short and simple proof. This should be compared with the (equally short) proof in \cite{loeb:05} that the contrapositive model existence theorem for decidable (and countable) classical propositional theories is equivalent to the Fan Theorem. (Note that the direction completeness $\Rightarrow$ Fan of Proposition \ref{theorem: fan equivalence} can be carried out in much weaker meta-theories  then IZF.
% in particular in the meta-theory of \cite{loeb:05}.)

\subsubsection{Chase-complete sets of models for regular theories} \label{Subsection: Regularbot theories}%%%%%%%%%%%%%%%%%%%%%%%%%%%%%%%%%%%%%%%%%%%%%%%%%

%$\chi$-closed $\chi$-complete $\operatorname{Ch}$-complete

Let $\Sigma$ be a single-sorted relational signature. That is, $\Sigma$ is a arbitrary set of relation symbols (with assigned arities), not assumed to be of a particular size nor discrete. 
%
%The assumption that it is relational is not essential; we prefer to avoid function symbols altogether, but it would make little difference if $\Sigma$ had constants, and we allow ourselves to extend signatures with constants from diagrams in the usual way.
%
%
%
%\cite{forsselllumsdaine} also shows that regular theories are, independently of size or discreteness constraints, ``faintly reflective'' in the following sense.
\begin{definition}\label{definition: chase functor}
Let \theory\ be a theory over $\Sigma$, and \cat{S} a class of $\Sigma$-diagrams. We say that  a functor  \funktor{\operatorname{Ch}}{\cat{S}}{\cat{S}}  is a \emph{chase functor} if for all diagrams $(D,F)\in \cat{S}$
\begin{enumerate}
\item $(D,F)\subseteq \operatorname{Ch}(D,F)$, naturally in $(D,F)$; 
\item  $\operatorname{Ch}(D,F)\vDash \theory$; and
\item
for any regular formula $\FIC{\alg{x}}{\psi}$ and  $\alg{d}\in D^{l(\alg{x})}$ such that $ \operatorname{Ch}(D,F)\vDash \psi[\alg{d}/\alg{x}]$, there exists a regular formula $\FIC{\alg{x}}{\phi}$ such that  $(D,F)\vDash \phi[\alg{d}/\alg{x}]$ and $\phi\vdash^{\theory}_{\alg{x}}\psi$.
\end{enumerate}
\end{definition}
%
%Strictly speaking this modifies the treatment  
%Let \theory\ be a regular theory over $\Sigma$.   Recall from Section \ref{subsection: models and diagrams} the assumption that the domains of diagrams are bounded relations to \thry{N}. With that, we extract the following from \cite{forsselllumsdaine}.

\begin{proposition}\label{theorem: chase completeness for regular theories} \label{proposition: chase completeness for regular theories}
Let \theory\ be a regular theory. Then there exists a chase functor \funktor{\operatorname{Ch}}{{\operatorname{Diag}_b(\Sigma)}}{\operatorname{Diag}_b(\Sigma)}. 
\begin{proof}
A straightforward modification of the construction in \cite{forsselllumsdaine}. (The modification being that,   working with  (bounded) diagrams rather than structures, we can have the components $c_{(D,F)}:(D,F)\to<125>\operatorname{Ch}(D,F)$ of the natural transformation $c:1_{\operatorname{Diag}_b(\Sigma)}\to \operatorname{Ch}$ be inclusions of diagrams, rather than homomorphisms of structures).
%
%   The slight modification here of requiring that there is an inclusion $(D,F)\subseteq \operatorname{Ch}(D,F)$ instead of just a homomorphism,  is made possible by the restriction to bounded diagrams.
\end{proof}
\end{proposition}
%
%We fix the construction of $(D\cup C, F')$ from \cite{forsselllumsdaine} and  write $\operatorname{Ch}(D,F)=(D\cup C, F')$. This construction of  is functorial in $(D\cup C, F')$ \cite{forsselllumsdaine}:
%
%\begin{proposition}\label{proposition: chase functoriality}
%Let \funktor{h}{(D,F)}{(C,E)} be a homomorphism of diagrams. Then there exists a homomorphism \funktor{\hat{h}}{\operatorname{Ch}(D,F)}{\operatorname{Ch}(C,E))} such that the square
%
%\[\bfig
%\square/>`<-`<-`>/<1000,400>[\operatorname{Ch}(D,F)`\operatorname{Ch}(C,E))`(D,F)`(C,E);\hat{h}`ch`ch`h]
%\efig\]
%commutes, where $ch$ are the homomorphisms induced by the inclusions. 
%\end{proof}
%\end{proposition}
%
 In general, we say that a collection of diagrams is \emph{closed under chase} if there is some chase functor \mo{Ch} under which it is closed. The construction of \cite{forsselllumsdaine} shows e.g.\ that if $\Sigma$ is a discrete signature,  then discrete diagrams are closed under chase.  In a classical meta-theory the stronger property holds that for a regular theory the category of models is weakly reflective in the category of structures (cf. e.g.\ \cite{adamekandrosicky:94}, \cite{makkai:90}).  We return to this and the case of enumerable\ signatures and theories in Section \ref{Subsection: Countable coherent theories and fan}. 
 %Finally, we remark in passing that if a coherent theory has a chase functor, then it is a regular theory. 

% construction of $\operatorname{Ch}(D,F)$ from $(D,F)$ shows that the category of models for a regular theory is weakly reflective in the category of structures (cf. e.g.\ \cite{adamekandrosicky:94}, \cite{makkai:90}). %That is to say, that any homomorphism from $(D,F)$ to a model diagram lifts to a (not necessarily unique) homomorphism from $operatorname{Ch}(D,F)$. 
%We display this for reference. 
%
%\begin{proposition}(AC)\label{proposition: chase weak refl} 
%
%^\noindent Let $(D,F)$ be a pre-diagram.   For any homomorphism \funktor{h}{(D,F)}{(D',F')} into a \theory-model diagram $(D',F')$ there exists a homomorphism $%\funktor{\hat{h}}{\operatorname{Ch}(D,F)}{(D',F') }$ such that 
%\[\bfig
%\ptriangle/>`<-`<-/<600,400>[\operatorname{Ch}(D,F)`(D',F')`(D,F);\hat{h}`ch`h]
%\efig\]
%commutes.
%
%\end{proposition}
%
%

For purposes of Section \ref{section: joyals theorem} we would, for given signature $\Sigma$ and regular theory \theory, like to restrict to a \emph{set} \cat{S} of diagrams  which is nevertheless ``rich enough'' for our purposes. Mainly, this involves being closed under chase, but we add some further conditions. Say that a diagram is %
%a \emph{fact extension} of a diagram $(D,F)$ if it is generated by a pre-diagram of the  form $(D, F\cup \pair{\syntob{\alg{x}}{\phi}, \alg{d}})$, where $\phi$ is a Horn formula and \alg{d} is a list of  elements of $D$. Say that a diagram is 
a \emph{finitary extension} of a diagram $(D,F)$ if it is generated by a pre-diagram of the  form $(D\cup \operatorname{Im}(\alg{c}), F\cup \pair{\syntob{\alg{x},\alg{y}}{\phi}, \alg{d},\alg{c}})$, where $\phi$ is a Horn formula, \alg{d} a (possibly empty) list of elements of $D$, and \alg{c} is a (possibly empty) list of  elements such that \Image{\alg{c}} is finite and disjoint from $D$. Say that a collection \cat{S} of diagrams is \emph{chase-complete}
%
%\footnote{For lack of a better term. No other meanings or uses of ``replete'' are intended.}
%
 if: the empty diagram is in \cat{S}; \cat{S} is closed under finitary extensions (up to isomorphism);    \cat{S} is closed under chase; and finally we add that, for any finite list of diagrams in  $ \cat{S}$ there exists mutually disjoint isomorphic copies in \cat{S} of those diagrams. We say that a collection \cat{M} of model diagrams is \emph{chase-complete} if \cat{M} is the collection of model diagrams in a chase-complete category of diagrams. By Proposition \ref{theorem: chase completeness for regular theories}, $\operatorname{MDiag}_b(\Sigma)$ is chase-complete. 
With reference to  the chase functor construction of \cite{forsselllumsdaine},  the restriction to a small and chase-complete subcategory of diagrams can be done e.g.\ by building a set \cat{U} based on the natural numbers and the syntax of the theory closed under finite lists; and then consider the bounded diagrams whose domains are subsets of $\cat{U}$. We leave the details. 
%we leave the details to the taste of the reader. 
For reference, we state, then:  
% 
%

%
%We fix the following terminology.
%\begin{definition}\label{definition: conservative, replete, etc}
%Let \theory\ be a regular (coherent) theory over a signature $\Sigma$. Let \cat{S} be a small, full subcategory of diagrams for $\Sigma$ and \cat{M} the full subcategory of \theory-model diagrams.  
%\begin{enumerate}
%\item \cat{M} is \emph{conservative} or \emph{complete} if  for every regular (coherent) sequent $\sigma$ over $\Sigma$, if $\sigma$ is true in all diagrams in \cat{M} then $\theory\vdash\sigma$.
%\item \cat{M} is \emph{geometrically conservative}  if  for every geometric sequent $\sigma$ over $\Sigma$, if $\sigma$ is true in all diagrams in \cat{M} then $\theory\vdash\sigma$.
%\item \cat{M} is \emph{strongly conservative} if it is conservative, and, moreover, for every $(D,F)\in \cat{M}$, the set of $\theory_{(D,F)}$-model diagrams the reducts of which are in \cat{M} is conservative for $\theory_{(D,F)}$.
%\item \cat{S} is \emph{replete} if: the empty diagram is in \cat{S}; for any finite list of diagrams in  $ \cat{S}$ there exists mutually disjoint isomorphic copies in \cat{S} of those diagrams; \cat{S} is closed under finite extensions (up to isomorphism); and   \cat{S} is closed under chase. We say that \cat{M} is \emph{replete} if it is the category of \theory-model diagrams in a replete category of structures. 
%\end{enumerate}
%\end{definition}
%
%We display for reference the conclusion of the current section.
%
\begin{theorem}\label{theorem: Reg theories have replete sets of models}\label{corollary: completeness for regular}
Every regular theory \theory\ has a chase-complete category \cat{M} of model diagrams. If the signature is discrete, then the theory has a chase-complete set of discrete model diagrams. 
%If the theory is enumerable, then it has a replete set of enumerable\ models.
%\begin{proof}
%With an eye on the construction of $\operatorname{Ch}(D,F)$ in \cite{forsselllumsdaine},   a replete set of diagrams for a signature $\Sigma$ can be constructed e.g.\ as follows. Form the set $U'$ by closing the (disjoint) sets of the signature, the variables, the logical symbols, and \thry{N} under finite lists. Let $U$  be the set of binary relations between $U'$ and \thry{N} with bounded range. And let \cat{M} be the set of diagrams $(D,F)$ with $D\in U$, tagged with the least upper bound of the domain.   We leave the details to the taste of the reader. 
 %Geometric conservativeness follows similarly by considering the pre-diagram $(\alg{s},\pair{\syntob{\alg{x}}{\phi}, \alg{s}})$.
%\end{proof}
\end{theorem}

We say that  \cat{M} is \emph{conservative} for a class \cat{K} of sequents if  for every sequent $\sigma$ in \cat{K}, if $\sigma$ is true in all diagrams in \cat{M} then $\theory\vdash\sigma$. (If \cat{K} is left implicit it is understood to be all sequents of the fragment of the theory). We say that \cat{M} is \emph{strongly conservative} for \cat{K} if it is conservative for \cat{K}, and, moreover, for every $(D,F)\in \cat{M}$, the set of $\theory_{(D,F)}$-model diagrams the reducts of which are in \cat{M} is conservative for $\theory_{(D,F)}$.

\begin{lemma}\label{lemma: regular strong completeness}
Let \theory\ be a regular theory and \cat{M} a chase-complete set of model diagrams. Then \cat{M} is strongly conservative.  
\begin{proof}
 Let \theory\ and $(D,F)$ be given, and let $\phi\vdash_{\alg{x}}\fins{\alg{y}}\psi$ be a normal form regular sequent over $\Sigma$ extended with $D$ as constants.  Assume this sequent to be true in all $\theory_{(D,F)}$-model diagrams the reducts of which are in \cat{M}. Replacing every occurrence of a constant from $D$ in $\phi$ with a fresh variable $z$, write $\phi=\phi'[\alg{d}/\alg{z}]$. Let \alg{s} be a list of fresh constants, disjoint from $D$, of the same length as \alg{x} such that $\operatorname{Im}(\alg{s})$ is finite. Then $\operatorname{Ch}(D\cup\operatorname{Im}(\alg{s}),F\cup\{\pair{\syntob{\alg{x},\alg{z}}{\phi'}, \alg{s}, \alg{d}}\})\vDash \theory_{D,F}$ and  $\operatorname{Ch}(D\cup\operatorname{Im}(\alg{s}),F\cup\{\pair{\syntob{\alg{x},\alg{z}}{\phi'}, \alg{s}, \alg{d}}\}\vDash \phi[\alg{s}/\alg{x}]$, so   $\operatorname{Ch}(D\cup\operatorname{Im}(\alg{s}),F\cup\{\pair{\syntob{\alg{x},\alg{z}}{\phi'}, \alg{s}, \alg{d}}\}\vDash \fins{\alg{y}}\psi[\alg{s}/\alg{x}]$. Whence $\theory_{(D\cup\operatorname{Im}(\alg{s}),F\cup\{\pair{\syntob{\alg{x},\alg{z}}{\phi'}, \alg{s}, \alg{d}}\}}$ proves the sequent $\top\vdash \fins{\alg{y}}\psi[\alg{s}/\alg{x}]$. Then $\theory_{(D,F)}$ proves the sequent $\phi[\alg{s}/\alg{x}]\vdash \fins{\alg{y}}\psi[\alg{s}/\alg{x}]$. With $\operatorname{Im}(\alg{s})$ being a finite set, we can conclude that $\theory_{(D,F)}$ proves the sequent $\phi\vdash_{\alg{x}} \fins{\alg{y}}\psi$.
\end{proof}
\end{lemma}
A similar argument shows also that a chase-complete \cat{M} is conservative for geometric sequents over the signature of \theory. We state this for reference.
\begin{lemma}\label{lemma: chasecomplete gives geometric conservative}
Let \theory\ be a regular theory and \cat{M} a chase-complete set of model diagrams. Then \cat{M} is conservative for geometric sequents over the signature of \theory.
\begin{proof}
See \cite{forsselllumsdaine}.
\end{proof}  
\end{lemma}

Finally, for the statement of Joyal's theorem in Section \ref{subsection: joyals theorem} we transfer the relevant results above to the case of structures for not necessarily purely relational signatures. This is a straightforward application of using the adjoint equivalence between diagrams and structures and of translating between signatures with function symbols and signatures without them, and we display it for reference. For theory \theory\ and model \alg{M}, the theory $\theory_{\alg{M}}$ of \alg{M} is defined as usual, so that $\modcat{\theory_{\alg{M}}}\simeq (\alg{M}\downarrow \modcat{\theory})$.
\begin{corollary}\label{corollary: strongly complete set of models}
Let \theory\ be a regular theory over an arbitrary signature $\Sigma$ (not necessarily purely relational). Then there exists a strongly complete set of models for \theory.  
\end{corollary}

\subsubsection{Enumerable coherent theories and Fan}\label{Subsection: Countable coherent theories and fan}

%Recall that by a \emph{enumerable\ theory} \theory\ we mean a theory such that the signature is a semi-decidable subset of a countable set and the set of axioms of \theory\ is semi-decidable. Similarly, by a \emph{enumerable\ model diagram} of a enumerable\ theory we mean a model diagram $(D,F)$ such that $D$ is a semi-decidable subset of a countable set and $F$ is semi-decidable (thus $\theory_{(D,F)}$ is enumerable\ in our sense). 
%
The construction of the functor \mo{Ch} of \cite{forsselllumsdaine} relied upon in the previous section involves applying all axioms of the theory simultaneously at each step.  (In that sense it could be said to be a ``simultaneous chase''.)  
In the enumerable setting one can, instead, apply a single axiom in each step. 
%allows for using a ``fair sequential chase'' instead. 
With  disjunctions  allowed in the axioms, this produces  a finitely branching  tree of structures, instead of a sequence of structures. Passing from regular to coherent theories, we therefore need the Fan theorem (with decidable bar\footnote{For statement and basic equivalents of the Fan theorem see e.g.\ \cite[Sect.\ 7]{troelstra:73}, where the relevant principle is named FAN$_\textnormal{D}$.}) to prove completeness. 
The construction in this case is akin to e.g.\ \cite{costelombardiroy:01} (and, as mentioned, the resulting proposition known), and we only outline it. Recall that by positive coherent we mean the coherent fragment without the logical constant $\bot$.   	 
\begin{proposition}[Fan]\label{lemma: completeness for sd coherent}\label{proposition: completeness for sd coherent}
Enumerable positive coherent theories are complete with respect to enumerable\ model diagrams.
\begin{proof}
Let \theory\ be a enumerable\ positive coherent theory over a relational signature $\Sigma$, assumed to be axiomatized by sequents on normal form.
%
%Notice that a Horn formula-in-context \FIC{\alg{z}}{\theta} can be seen as specifying a particular finite diagram, namely the least diagram for which the domain consists of the variables \alg{z} and $\theta$ is true. 
%
%Recall that we assume, without loss of generality, that the axioms of a positive coherent theory are written as sequents where the antecedent is a Horn formula and the consequent is a finite, inhabited disjunction of regular formulas. 
%
An \emph{application} of such an axiom
\[\theta\vdash_{\alg{x}}\bigvee_{1\leq i\leq n} \fins{\alg{y_i}}\psi_i\]
to a diagram $(D,F)$   is a function $f:\alg{x}\rightarrow D$ such that $\pair{\syntob{\alg{x}}{\theta},f(\alg{x})}\in F$. Such an application induces $n$ children $(D_i,F_i)$, where   $(D_i,F_i)$ is the least extension of $(D,F)$ containing a list $\alg{y'}_i$ of distinct fresh elements of the same length as $\alg{y}_i$ and such that $\psi_i[f(\alg{x})/\alg{x},\alg{y'}_i/\alg{y}_i]$ is true. 

Given a finite diagram $(D,F)$, build a finitely branching tree \cat{T} of finite diagrams with root $(D,F)$ e.g.\ as follows. With \theory\ enumerable\ we can find finite subtheories  $\theory_n$ such that $\theory_0\subseteq \ldots \theory_n \subseteq \ldots \bigcup_{n\in \thry{N}}\theory_n=\theory$. Build a sequence of finite trees $\cat{T}_n$ with $\cat{T}_0$ consisting just of $(D,F)$ and $\cat{T}_n$ an initial subtree of $\cat{T}_{n+1}$ as follows.   For each leaf $(D',F')$ of the tree $\cat{T}_n$, list all possible applications $a_1,\ldots,a_m$ of $\theory_n$ to that leaf. Apply $a_1$ to $(D',F')$. This induces a finite number of children which are extensions of $(D',F')$. Apply $a_2$ to those, $a_3$ to the children induced by that again,  and proceed until the list of applications runs out. This produces a finite tree $\cat{T'}_{(D',F')}$ with $(D',F')$ as its root. Then $\cat{T}_{n+1}$ is obtained by appending  $\cat{T'}_{(D',F')}$ to each leaf $(D',F')$ of $\cat{T}_n$. Let \cat{T} be the union of the $\cat{T}_n$.   Notice that: 1) $\cat{T}_n$ constitutes (or translates to) a dynamical cover (c.f.\ \cite{costelombardiroy:01}) of $(D,F)$ with respect to $\theory$; and 2) the union of the diagrams along a path of \cat{T} is a (enumerable) \theory-model diagram.

Now, let \[\phi\vdash_{\alg{x}}\bigvee_{ i\leq n} \fins{\alg{y_i}}\psi_i\] be a positive coherent sequent (on normal form, without loss of generality). Assume it is true in all \theory-model diagrams. Let $(D,F)$ be the finite diagram presented by the formula-in-context \FIC{\alg{x}}{\phi}; that is the least diagram for which the domain consists of a finite set bijective with $\mo{Im}(\alg{x})$, we can write $\bar{\alg{x}}$, and $\phi[\bar{\alg{x}}/\alg{x}]$ is true. Construct the tree \cat{T} on $(D,F)$ as above. Define the subset of nodes
\[B=\cterm{(D',F')}{\fins{i\leq n}\fins{\alg{d}\in D'}\psi_i[\bar{\alg{x}}/\alg{x},\alg{d}/\alg{y}_i]\in F'}\]
This is a decidable subset, and since the sequent is true in all models, and thus in all paths, it is a bar. Thus by the Fan theorem, it is a universal bar. Accordingly,  there is an $n$ such that every leaf node in $\cat{T}_n$ is in $B$. And since $\cat{T}_n$ is a dynamical cover, the sequent is provable  in $\theory$.
\end{proof}
\end{proposition}

The fallible Kripke semantics of e.g.\ \cite{veldman:76} allows for \emph{exploding} nodes, in the form of nodes that force $\bot$. Such nodes  force all other formulas as well. Similarly we could define an exploding structure as one that interprets $\bot$ as true. We prefer to look at this through the lense of Morleyization; a structure for $\theory^m$ yields a structure for \theory\ in which some of the logical constants are interpreted non-standardly. In particular, if \theory\ is a coherent theory and $\theory^m$  its positive coherent Morleyization, then a $\theory^m$-model diagram $(D,F)$ induces an interpretation of the formulas of \theory, by letting the extension of a formula $\phi$ be the extension of the corresponding predicate $P_{\phi}$. We say that this interpretation is \emph{exploding} if $P_{\bot}\in F$. And we refer to the interpretation of \theory\ in terms of the models of $\theory^m$ as possibly exploding or \emph{fallible} Tarski semantics. We will consider further modifications in the next section.    
%
%Informally, we say that a structure is exploding if it gives a non-standard interpretation of $\bot$. Formally, and for current purposes, we define the category of \emph{possibly exploding} model diagrams of a coherent (regular$_{\bot}$ theory \theory\ to be the category of model diagrams for its positive coherent (regular) Morleyization. We refer to the intepretation of a coherent (regular$_{\bot}$ theory with respect to such models as \emph{modified Tarski semantics}. 
%
We now have the following corollary of Proposition \ref{lemma: completeness for sd coherent}.
\begin{corollary}[Fan]\label{corollary: Coherent modified tarski completeness}
Enumerable coherent theories are complete with respect to fallible Tarski semantics. 
\end{corollary}
%
%We note also, as a scholium,  the enumerable\ counterpart of Proposition \ref{proposition: Reg theories have replete sets of models}. It does not require the Fan Theorem. 
%
%\begin{scholium}\label{proposition: SD Reg theories have replete sets of sd models}
%Every enumerable\ regular theory has a replete set of enumerable\ models.
%\begin{proof}[Sketch]
%
%\marginpar{[TODO: check, expand sketch?]}
%Let $(D,F)$ be a enumerable\ diagram.  In this case, the construction of $\operatorname{Ch}(D,F)$ can be done, in a matter of %speaking, by using a ``fair chase'' instead of a ``simultaneous chase''. That is, in outline, as follows. Generate in stages finite %sets of axioms of \theory\ and finite sets of facts of $(D,F)$. At each stage make all possible applications of the axioms obtained %to the facts obtained, and add the results to the facts obtained at the next stage. The model diagram obtained at the limit is %then enumerable   
%\end{proof}
%\end{scholium}
%
%
%[We will\marginpar{[TODO: align this with the relavant application and write out accordingly]} assume (in Section \ref{subsection: modified completeness} in particular) that a replete set of enumerable\ models for a enumerable\ regular theory consists of models diagrams the domains of which are left-bounded semi-decidable subsets of \thry{N} and otherwise unifromly organized in a suitable way.]
%
%
If \theory\ is an enumerable classical first-order theory, then its coherent Morleyization $\theory^m$ is an enumerable\ coherent theory. A model for $\theory^m$ can be regarded as an ordinary Tarski model for \theory\ satisfying  LEM; that is to say, every set that is the extension of a first-order formula of the language of the theory \theory\ must be complemented. Thus regarding a classical first-order theory as a first-order theory containing the LEM axiom scheme, we have as a consequence: 
\begin{corollary}[Fan]\label{corollary: completeness of classical fo theories}
Classical enumerable\ first-order theories are complete with respect to fallible Tarski semantics.
\end{corollary}
%
%
%\begin{remark}\label{Remark: on the completeness of cfol}
%Intuitionistic completeness theorems for classical logic have been shown (see, e.g., \cite{krivine:96}) working with a %notion of model which basically amounts to modified Tarski semantics, but in which the definition of validity is slightly %changed, and as a result the standard meaning of disjunction is lost. Keeping the classical meaning of the logical %connectives, the algorithmic content of completeness theorems for classical logic has been exposed, e.g., in \cite{kleene}, %but the metatheory is classical and in fact makes a strong use of K\"onig's lemma. Corollary \ref{corollary: completeness of %classical fo theories} shows, however, that a completely intuitionistic treatment is possible retaining the meaning of the %connectives, if one also allows exploding models.
%\end{remark}
%
\begin{remark}\label{Remark: on the completeness of cfol}
Corollary \ref{corollary: completeness of classical fo theories} provides an intuitionistic completeness theorem for classical logic provided the notion of model is relaxed to allow exploding models. Such theorems have been derived before, notably by Krivine (see \cite{krivine:96} and also \cite{berardiandvalentini:04}).
% However, these use a notion of semantics in Tarski structures that is not standard for all logical connectives. 
%
In fact, Krivine proves the model existence theorem for consistent classical first-order theories (which is intuitionistically stronger than the completeness theorem),  %which asserts that validity (in all structures, including possibly the exploding ones) implies provability. 
%Given a 
%consistent theory, he proves that there is a model 
%
%in which the standard semantics of all connectives, except disjunction, is preserved.
%
but with respect to models in which disjunction is non-standardly interpreted. In comparison, Corollary \ref{corollary: completeness of classical fo theories} retains the standard semantics for all connectives except $\bot$.
%
% but at the cost of proving the (weaker)  completeness theorem instead of the model existence theorem.
\end{remark}
Finally, we show that the use of the Fan theorem in Proposition \ref{lemma: completeness for sd coherent} is essential, and conclude:
\begin{theorem}\label{theorem: fan equivalence}
The completeness of enumerable\ positive coherent theories with respect to enumerable\ model diagrams is equivalent to the Fan theorem.
\begin{proof}
Let $F$ be a fan containing a decidable bar $\alg{B}$. We follow the notation in \cite[4.1]{troelstravandalen:88}.  Consider the theory $\theory$ over the signature consisting of a propositional variable $P_{\alg{n}}$ for each element $\alg{n} \in F$ and a propositional variable $B$, and whose axioms are:
\begin{enumerate}
\item $\top \vdash P_{\langle \rangle}$, where $\langle \rangle$ is the root of the Fan;
\item $P_{\alg{n}} \vdash \bigvee_{\alg{n}\ast m \in F} P_{\alg{n}\ast m}$;
\item $P_{\alg{n}} \vdash B$ for each $\alg{n} \in \alg{B}$.
\end{enumerate}
For each branch $\alpha$ in $F$ let $S_{\alpha}$ be the set of sentences  
\[S_{\alpha}=\cterm{P_{\alg{n}}}{\fins{x}\bar{\alpha}x=\alg{n}}\]
%\[S_{p}=\cterm{P(\alpha)}{\alpha\in p}\bigcup \cterm{B}{\fins{n}p(n)\in \alg{B}}\]
%
%Since $\alg{B}$ is a bar, $B\in S_p$. 
Let $\alg{M}=(D,F)$ be a enumerable\ \theory -model. Since it is enumerable\ we can find a branch $\alpha$ such that $S_{\alpha}\subseteq F$.
Since \alg{B} is a bar, we then have that $B\in F$. Thus $\alg{M}\vDash (\top\vdash B)$. By completeness, there is a proof of $\top\vdash B$ in \theory, with finitely many axioms, whence \alg{B} must be uniform. 
\end{proof}
\end{theorem}

Since  classical first-order logic is conservative over coherent logic (see e.g.\ \cite{negri:03}), we could have added the fallible Tarski-completeness of classical FOL as a third equivalent statement in  \ref{theorem: fan equivalence}. We proceed now to the theorem of Joyal by which one can add the fallible Kripke completeness of FOL as a fourth.
%
%
%From the proof of \ref{theorem: fan equivalence} and from Lemma \ref{3} we obtain:
%
%\begin{corollary}\label{corollary: fan equivalence disjunctive}
%Completeness with respect to model diagrams of s.d\ positive coherent theories is equivalent to the Fan theorem.
%\end{corollary}
%
%\begin{remark}
%Theorem \ref{1} (together with the results of \cite{mccarty:08}) and Corollary \ref{corollary: fan equivalence disjunctive} %display the need for exploding models (or Markov's principle) for fragments containing $\bot$ and the need for the Fan %theorem for fragments containing disjunction.
%\end{remark}
%
%[TODO: check that all these fragments make sense. expand remark]
%
%
Before doing so, however, we note, for use in  Section \ref{section: sheaf completeness},  that if \theory\ is a regular theory then the construction of Proposition \ref{lemma: completeness for sd coherent} yields a sequence of diagrams, the union of which is a model of \theory. That this construction, extended to general enumerable\ diagrams,  can be used to show that enumerable\ \theory\ models are weakly reflective in enumerable\ diagrams, as in the classical case, is rather expected and straightforward.
We therefore state the following for reference and without proof.    
%
%Note that if \theory\ is in fact a regular theory, then \cat{T} is a sequence of diagrams. And it is straightforward to see that a homomorphism from the root to a enumerable\ diagram can be lifted along this sequence. The construction can also be extended from a finite starting diagram $(D,F)$ to an enumerable\ one by considering increasingly large finite subdiagrams. Thus enumerable\ regular theories are not only replete with respect to enumerable\ model-diagrams, but Proposition \ref{proposition: chase weak refl} holds. We display this for reference. (The proof holds no surprises and is omitted.)  
%
\begin{proposition}\label{proposition: sd regular chase} 

\noindent Let \theory\ be an enumerable\ regular theory. There exists a chase-complete set \cat{S} of enumerable\ diagrams with a chase functor $\funktor{\operatorname{Ch}}{\cat{S}}{\cat{M}}$---where \cat{M} is the subcategory of model diagrams---which is moreover a weak reflection. That is,   
for  any homomorphism \funktor{h}{(D,F)}{(D',F')} in \cat{S} where $(D',F')\in \cat{M}$ there exists a homomorphism $\funktor{\hat{h}}{\operatorname{Ch}(D,F)}{(D',F') }$ such that 
\[\bfig
\ptriangle/>`<-`<-/<600,400>[\operatorname{Ch}(D,F)`(D',F')`(D,F);\hat{h}`c_{(D,F)}`h]
\efig\]
commutes.
\end{proposition}

\section{Joyal's theorem}\label{section: joyals theorem}\label{subsection: joyals theorem}%%%%%%%%%%%%%%%%%%%%%%%%%%%%%%%%%%%%%%%%%%%%%%%%%%%%%
%%%%%%%%%%%%%%%%%%%%%%%%%%%%%%%%%%%%%%%%%%%%%%%%%%%%%%%%%%%%%%%%%%%%%%%%%%%%%%%%%%%%%%%%%%%%%%%%%%%%%%%%%%%%%

%\subsection{Joyal's theorem}\label{subsection: joyals theorem}%%%%%%%%%%%%%%%%%%%%%%%%%%%%%%%%%%%%%%%%%%%%%%%%%%%%%
%%%%%%%%%%%%%%%%%%%%%%%%%%%%%%%%%%%%%%%%%%%%%%%%%%%%%%%%%%%%%%%%%%%%%%%%%%%%%%%%%%%%%%%%%%%%%%%%%%%%%%%%%%%%%

Let \theory\ be a coherent (or regular) theory. A coherent (regular) formula-in-context \syntob{\alg{x}}{\phi} induces an evaluation functor
\[\funktor{Ev_{\syntob{\alg{x}}{\phi}}}{\modcat{\theory}}{\Sets}\]
by $\alg{M}\mapsto \csem{\alg{x}}{\phi}^{\alg{M}}$. 
Mapping a formula-in-context to its corresponding evaluation functor defines (by soundness) a functor $\operatorname{Ev}:\synt{C}{\theory}\to<125>\Sets^{\modcat{\theory}}$, which we also call the evaluation functor, trusting that context will prevent confusion.
Since the coherent structure in a presheaf category is computed pointwise, the following is immediate and stated only for emphasis and reference.
\begin{lemma}\label{lemma: ev functor is coherent}
Let \theory\ be a coherent \textup{(}regular\textup{)} theory and \synt{C}{\theory} its coherent \textup{(}regular\textup{)} syntactic category. The functor
\[\operatorname{Ev}\ :\ \synt{C}{\theory}\to \Sets^{\modcat{\theory}}\]
which sends a formula to its corresponding evaluation functor is coherent \textup{(}regular\textup{)}.
\end{lemma}
We can now give the ``constructive content'' of the Kripke completeness theorem of A.\ Joyal---cf.\  \cite[Thm 6.3.5]{makkaireyes}---first in the form of the following theorem for  regular theories. Since the purpose is to give a constructive restatement of this classical theorem, we state it first in terms of arbitrary signatures and ordinary Tarski models. The proof is not in  essence dissimilar from the one in \cite{makkaireyes}.  Recall that by Corollary \ref{corollary: strongly complete set of models} there are strongly complete sets of models for regular theories. 

\begin{theorem}%[Joyal]
\label{theorem: joyals theorem}
Let $\Sigma$ be a single sorted theory, not restricted in size, nor necessarily discrete  (and possibly containing function symbols).  Let \theory\ be a regular theory over $\Sigma$, and let $\cat{M}$ be a full subcategory of $\modcat{\theory}$ such that \cat{M} is strongly conservative.
%\begin{enumerate}[i)]
%\item \theory\ is complete with respect to  $\cat{M}$, and 
%\item  for every $\alg{M}\in \cat{M}$ the theory $\theory_{\alg{M}}$ is complete with respect to model diagrams \textup{(}the reducts of which are\textup{)} in \cat{M}.
%\end{enumerate}
Then the functor
\[\operatorname{Ev}\ :\ \synt{C}{\theory}\to \Sets^{\cat{M}}\]
is a\textup{)} conservative and b\textup{)} whenever the pullback functor $\funktor{f^*}{\mathrm{Sub}_{\synt{C}{\theory}}(B)}{\mathrm{Sub}_{\synt{C}{\theory}}(A)}$ induced by a morphism $f\! :\! A\to<125>B$ in \synt{C}{\theory} has a right adjoint $\forall_f$ we have   for all $S\in\sublat{\synt{C}{\theory}}{A}$ that $Ev(\forall_f(S))=\forall_{Ev(f)}(Ev(S))$.
\begin{proof} a) For formulas \syntob{\alg{x}}{\phi} and \syntob{\alg{x}}{\psi}, if  $\operatorname{Ev_{\syntob{\alg{x}}{\phi}}}(\alg{M})\subseteq \operatorname{Ev_{\syntob{\alg{x}}{\psi}}}(\alg{M})$ for all $\alg{M}\in\cat{M}$ then $\alg{M}\vDash (\phi\vdash_{\alg{x}}\psi)$ for all  $\alg{M}\in\cat{M}$, whence  $\phi\vdash^{\theory}_{\alg{x}}\psi$ by completeness.

b) The non-trivial direction is $Ev(\forall_f(S))\supseteq\forall_{Ev(f)}(Ev(S))$. It suffices to consider a situation 
\[\bfig
\square/`^{(}->`^{(}->`>/<2000,300>[S=\syntob{\alg{x}}{\theta}`\forall_f(S)=\syntob{\alg{y}}{\gamma}`A=\syntob{\alg{x}}{\phi}`B=\syntob{\alg{y}}{\psi};```f=|\syntob{\alg{x},\alg{y}}{\lambda}|]

\efig\]
in \synt{C}{\theory}, where $\theta\vdash^{\theory}_{\alg{x}}\phi$ and $\gamma\vdash^{\theory}_{\alg{y}}\psi$.
Applying the functor Ev and evaluating at a model $\alg{M}$ we have 
\[\bfig
\square/`^{(}->`^{(}->`>/<2000,300>[\cterm{\alg{d}}{\alg{M}\vDash\theta(\alg{d})}`\cterm{\alg{c}}{\alg{M}\vDash\gamma(\alg{c})}`\cterm{\alg{d}}{\alg{M}\vDash\phi(\alg{d})}`\cterm{\alg{c}}{\alg{M}\vDash\psi(\alg{c})};```\operatorname{Ev}_f(\alg{M})=\cterm{\alg{d},\alg{c}}{\alg{M}\vDash\lambda(\alg{d},\alg{c})}]

\efig\]
%
%
%\[\bfig
%\square/`>`>`>/<2500,500>[\evalmod{S}{D,F}=\cterm{[\alg{d}]}{(D,F)\vDash\theta(\alg{d})}`\evalmod{\forall_f(S)}{D,F}=\cterm{[\alg{c}]}{(D,F)\vDash\gamma(\alg{c})}`\evalmod{A}{D,F}=\cterm{[\alg{d}]}{(D,F)\vDash\phi(\alg{d})}`\evalmod{B}{D,F}=\cterm{[\alg{c}]}{(D,F)\vDash\psi(\alg{c})};```\evalmod{f}{D,F}=\cterm{[\alg{d}],[\alg{c}]}{(D,F)\vDash\lambda(\alg{d},\alg{c})}]
%
%
%\efig\]
%
Let $\alg{c}\in \forall_{\operatorname{Ev}_{f}}(\operatorname{Ev}_{\syntob{\alg{x}}{\theta}})(\alg{M})\subseteq \operatorname{Ev}_{\syntob{\alg{y}}{\psi}}(\alg{M})$. Accordingly, for all $g\! :\! \alg{M}\to<125>\alg{N}$  in $\cat{M}$ we have:
\begin{equation}\label{Eq: Joyals-forall eq}\evalmod{f}{(\alg{N})}^{-1}(\operatorname{Ev}_{\syntob{\alg{y}}{\psi}}(g)(\alg{c}))\subseteq \operatorname{Ev}_{\syntob{\alg{x}}{\theta}}(\alg{N}).\end{equation}
We show $\alg{M}\vDash\gamma(\alg{c})$. Let $g\! :\! \alg{M}\to<125>\alg{N}$ be a morphism in $\cat{M}$, with  $\alg{N}'$ the corresponding $\theory_{\alg{M}}$-model.
By (\ref{Eq: Joyals-forall eq}) we have
\begin{equation}\label{Eq: Joyals-forall seq}
\alg{N}'\vDash (\lambda[\alg{c}/\alg{y}]\vdash_{\alg{x}}\theta)
\end{equation}
Thus, the sequent (\ref{Eq: Joyals-forall seq}) is true in all $\theory_{\alg{M}}$-models corresponding to homomorphisms from $\alg{M}$ in \cat{M},
%
%objects in $\left((D,F)\downarrow\cat{M}\right)$,
%
and therefore provable in $\theory_{\alg{M}}$, by the assumption of strong completeness. By Lemma \ref{lemma: getting rid of constants from the model}, there is a regular formula $\xi$ in context \alg{y} such that $\alg{M}\vDash\xi(\alg{c})$ and $\theory$ proves the sequent $(\xi\wedge\lambda\vdash_{\alg{x},\alg{y}}\theta)$. But then, since 
%$S=\syntob{\alg{x}}{\theta}$ and 
$\forall_f(\syntob{\alg{x}}{\theta})=\syntob{\alg{y}}{\gamma}$, we have that \theory\ proves the sequent $(\xi\wedge \psi\vdash_{\alg{y}}\gamma)$. Whence $\alg{M}\vDash \gamma(\alg{c})$.
\end{proof}
\end{theorem}
%
%The formulation of Theorem \ref{theorem: joyals theorem} in terms of coherent theories and models (and hence the classical statement of Joyal's theorem on p.\pageref{Theorem-non: Joyal}) is a straightforward corollary (see \ref{corollary: classical joyals} below). 
%

It is convenient to have a name for the property proved in Theorem \ref{theorem: joyals theorem}. Following e.g.\ \cite{butz:02} (at least for the first notion):
\begin{definition}\label{definition: subheyting}
We say that a functor $\funktor{F}{\cat{C}}{\cat{D}}$ from a coherent category to a Heyting category is \emph{conditionally Heyting} if it is coherent and preserves any right adjoints to pullback functors that might exist in \cat{C}. If $F$ and \cat{C} are regular, we say $F$ is \emph{conditionally sub-Heyting}  if it is regular and preserves any right adjoints to pullback functors that might exist in \cat{C}.
\end{definition}

The equivalence $\modcat{\theory}\simeq \moddiag{\theory}$ induces an equivalence $\Sets^{\modcat{\theory}} \simeq  \Sets^{\moddiag{\theory}} $. Accordingly,  if \cat{M} is a strongly complete set of model diagrams, the composite 
\[  \synt{C}{\theory}\to^{\operatorname{Ev}} \Sets^{q(\cat{M})} \simeq  \Sets^{\cat{M}}\]  
is conservative and  conditionally sub-Heyting. When returning to working with diagrams in the sequel, we shall consider this functor, also under the name $\operatorname{Ev}$.  

 As a corollary Theorem \ref{theorem: joyals theorem} we have the following version for coherent theories.  Strongly complete sets of models do not in general exist for coherent theories. By Proposition\ref{lemma: completeness for sd coherent}, however, they do for enumerable\ positive coherent theories under the assumption of  the Fan theorem.

%(The notion of strongly conservative set or category of models for a coherent theory is the same as in Definition \ref{definition: conservative, replete, etc} replacing %``regular'' by ``coherent''.)
%
\begin{corollary}\label{corollary: classical joyals}
Let \theory\ be a coherent theory, with \synt{C}{\theory} its coherent syntactic category, and suppose  \cat{M} is a strongly conservative category of \theory-model  diagrams. 
%satisfying conditions (i) and (ii) of Theorem \ref{theorem: joyals theorem}. 
%
Then the evaluation functor 
\[\operatorname{Ev}\ :\ \synt{C}{\theory}\to \Sets^{\cat{M}}\]
is a conservative and conditionally Heyting functor.
% which preserves right adjoints to pullback functors whenever they exist. 
\begin{proof}The proof of \ref{theorem: joyals theorem} can be repeated for this case. Alternatively,  consider the regular Morleyization  $\theory^m$ of the coherent theory \theory.  Notice, first, that \cat{M} can then be considered as a full subcategory of $\operatorname{Mod}(\theory^m)$, and as such it is then strongly conservative for $\theory^m$. Then notice that evaluation restricted to \cat{M} yields a coherent functor from $\synt{C}{\theory}\simeq \synt{C}{\theory^m}$ to $\Sets^{\cat{M}}$.
\end{proof}
\end{corollary}

%\subsection{Morleyization and modified completeness}\label{section: morleyization and modified completeness}%%%%%%%%%%%%
%%%%%%%%%%%%%%%%%%%%%%%%%%%%%%%%%%%%%%%%%%%%%%%%%%%%%%%%%%%%%%%%%%%%%%%%%%%%%%%%%%%%%%%%%%%%%%%%%%%%%%%%%%%%%%%%%%%%%%%%

\section{Sheaf completeness}\label{section: sheaf completeness}

\subsection{Modified completeness}\label{subsection: modified completeness}%%%%%%%%%%%%%%%%%%%%%%%%%%%%%%%%%
%%%%%%%%%%%%%%%%%%%%%%%%%%%%%%%%%%%%%%%%%%%%%%%%%%%%%%%%%%%%%%%%%%%%%%%%%%%%%%%%%%%%%%%%%%%%%%%%%%%%%%%%%%%%%%%%%%

Loosely and informally, let us say that a model is \emph{modified} if some connectives are interpreted in a non-standard way, and \emph{standard} otherwise. We say that it is \emph{fallible} if the only connective treated non-standardly is $\bot$.  
%
%
% It is well-known, or at least ``folklore'', that regular logic is (constructively) complete with respect to standard Tarski models. (For conceptual reasons and  for presentational completeness, we present a short proof of this in Section \ref{section: modified completeness}). 
%
%Thus from Theorem \ref{theorem: joyals theorem} we have the following.
%
%\begin{corollary}\label{corollary: disjunction and falsum free joyals}
%Let \theory\ be a theory in the fragment consisting of $\top$, $\wedge$,$\exists$,$\rightarrow$, and $\forall$, and  
%$\theory^m$ its regular Morleyization. Then 
%\[\operatorname{Ev}\ :\ \synt{C}{\theory}\simeq \synt{C}{\theory^m}\to \Sets^{\operatorname{MDiag}(\theory^m)}\]
% yields a conservative model of \theory\ in the presheaf category $\Sets^{\operatorname{MDiag}(\theory^m)}$. 
%\end{corollary}
%
%We refer to $\theory^m$ as the ``Morleyization'' of \theory\ (cf. \cite{elephant1}). Since this axioms specify (up to %equivalence of categories) the internal coherent (regular$_{\bot}$) theory of \synt{C}{\theory} (cf.\ \cite{elephant1}) we %have that
%
%\[\synt{C}{\theory}\simeq \synt{C}{\theory^m}\]
%
%This amounts (in familiar ways) to Kripke completeness for that fragment, as we shall briefly return to in Section %\ref{section: veldman etc thms}.
%
%Thus for a regular theory \theory\ the category \moddiag{\theory} satisfies the conditions (i) and (ii) of Theorem \ref{theorem: joyals theorem}. 
%
%In particular, fallible Tarski models are modified in treating $\bot$ non-standardly. 
%
It is a corollary of Theorem \ref{corollary: completeness for regular} and Lemma \ref{lemma: regular strong completeness} that regular$_{\bot}$ theories are complete with respect to  fallible Tarski semantics.
%
%As mentioned in Section \ref{subsubsection: morleyization}, possibly exploding Tarski models of a regular$_{\bot}$ theory \theory\ are the same thing, or can be seen as,  ordinary models of its regular Morleyization $\theory^m$. Therefore, it follows from the completeness of regular logic  that regular$_{\bot}$ logic is complete with respect to  Tarski semantics modified in  terms of allowing  exploding models.  
%
%
Now, if \theory\ is, say, a regular$_{\bot}$ theory,  $\theory^m$ it's regular Morleyization, and \cat{M} is a small, full, and strongly conservative subcategory of $\operatorname{MDiag}(\theory^m)$,  we have that the evaluation functor 
\[\operatorname{Ev}\ :\ \synt{C}{\theory}\simeq \synt{C}{\theory^m}\to \Sets^{\cat{M} }\]
is regular (in particular). But it does not preserve the initial object, as $\evalmod{\syntob{}{P_{\bot}}}{-}$ is not the constant empty functor $0$. Thus it can be viewed as a conservative fallible presheaf model of \theory. We shall obtain a conservative standard sheaf model by taking sheaves with respect to the least coverage (on \op{\cat{M}}) so that $\evalmod{\syntob{}{P_{\bot}}}{-}$ is identified with $0$. Accordingly,  we obtain a model of \theory\ in a closed subtopos (in the sense of \cite[A4.5.3]{elephant1}) of $\Sets^{\cat{M}}$. Similarly, if \theory\ is a coherent theory and $\theory^m$ its regular Morleyization, $\operatorname{Ev}\ :\ \synt{C}{\theory}\to \Sets^{\cat{M} }$ does not preserve finite disjunctions. A conservative standard model will be obtained by sheafifying with respect to the least coverage such that finite disjunctions  are preserved. A  conservative model will also be given by a slightly stronger coverage given (in part) in terms of binary trees and which is, in that sense, akin to a (fallible) Beth model. Classically, or in a enumerable setting, the latter two coverages are equivalent. In a countable setting, they also give rise to a Beth-completeness theorem, establishing a link between the least coverage forcing a standard interpretation and Beth semantics.   

The coverages are given in terms of sieves on \op{\cat{M}}, and thus cosieves on \cat{M}.
%
% Since this is rather transparent, however, we shall speak of sieves and coverages on \cat{M}, although we strictly speaking mean on  \op{\cat{M}}.
%
%Although we shall speak of \cat{M}  rather than \op{\cat{M} }, we avoid speaking of cosieves or ``cocoverage'' and trust that no confusion arises from this.    
Explicitly, then, let \theory\ be a theory in a fragment with $\bot$,  $\theory^m$ its regular Morleyization, and \cat{M} a small full subcategory of \moddiag{\theory^m}.
 Let the \emph{exploding coverage} ${E}$ be the coverage which assigns to each $(D,F)\in \cat{M}$ the set of cosieves $E(D,F)=\cterm{\emptyset}{P_{\bot}\in F}$.  
This is a coverage (in the sense of \cite[A2.1.9, C2.1.1]{elephant1}) since if $S\in E(D_1,F_1)$ and  $f:(D_1,F_1)\to<150> (D_2,F_2)$ is a homomorphism, then $S=\emptyset$ and $P_{\bot}\in F_1$, whence $P_{\bot}\in F_2$, since it is preserved by $f$,  so $\emptyset\in E(D_2,F_2)$.  
We then have the following addition to Theorem \ref{theorem: joyals theorem}.

\begin{proposition}\label{proposition: modified joyals theorem for regular}
Let \theory\ be an \textup{(}at least\textup{)} regular$_{\bot}$ theory, $\theory^m$ its regular Morleyization, \cat{M} be a strongly conservative, small, full subcategory of $\operatorname{MDiag}(\theory^m)$  and $E$ be the exploding coverage. Then the evaluation functor factors through sheaves
\[\bfig
\qtriangle/>`>`<-^{)}/<500,300>[\synt{C}{\theory}\simeq\synt{C}{\theory^m}`\Sets^{\cat{M}}`\Sh{\op{\cat{M}},E};\operatorname{Ev}``]
\efig\]
So that $\operatorname{Ev}:\synt{C}{\theory}\to<150>\Sh{\op{\cat{M}},E}$ is conservative, conditionally sub-Heyting, and preserves the initial object.
\begin{proof} By Theorem \ref{theorem: joyals theorem} it remains only to show that $\operatorname{Ev}$ factors through \Sh{\op{\cat{M}},E} and that $\evalmod{\syntob{}{{\bot}}}{-}$ is terminal in \Sh{\op{\cat{M}},E}. First, $\evalmod{\syntob{\alg{x}}{\phi}}{-}$ is a sheaf since if
 $S \in E(D,F)$, then $S$ is empty and $(D,F)$ is exploding, which means that $\evalmod{\syntob{\alg{x}}{\phi}}{D,F}=\{\ast\}$. 
%It is straightforward to see that the functor $ev:\synt{C}{\theory}\to<150>\op{\moddiag{\theory}}$ %is regular and conservative, the latter by Theorem \ref{1}. We show that the inclusion %$\topo{E}\hookrightarrow \Sets^{\moddiag{\theory}}$ preserves covers. Let $\Phi:H\to<150>G$ be a %natural transformation of sheaves, with $R\subseteq G$ the presheaf image of $\Phi$. Then the %closure of $R$ is
%
%\[\overline{R}(D,F)=\cterm{x\in G(F,D)}{\fins{S\in E(F,D)}\alle{f\in S}G(f)(x)\in R(cod(f))}\]
%
%Let $x\in \overline{R}(D,F)$. Then there exists   $S\in E(F,D)$ such that for all $f\in S$ we have %$G(f)(x)\in R(cod(f)$. If $S$ is maximal then $x\in R(D,F)$. If $S$ is empty then $P_{\bot}\in F$ whence %$H(D,F)=G(D,F)=1$ so $x\in R(D,F)$. Thus $\overline{R}(D,F)=R(D,F)$.
%
Second, $\evalmod{\syntob{}{{\bot}}}{D,F}=\cterm{\ast}{P_{\bot}\in F}$ which is the initial sheaf in \Sh{\op{\cat{M}},E}.
\end{proof}
\end{proposition}
Note that classically the standard (i.e.\ non-exploding) models in \cat{M} are dense (in the sense of \cite{elephant1}), so that then $\Sh{\op{\cat{M}},E}\simeq \Sets^{\cat{M}^s}$ where $\cat{M}^s$ is the full subcategory of standard models.

Next, let \theory\ be a theory in a fragment with $\vee$ and $\bot$, with $\theory^m$ its regular Morleyization and \cat{M} a small full subcategory of \moddiag{\theory^m}. Again, 
\[\operatorname{Ev}\ :\ \synt{C}{\theory}\simeq\synt{C}{\theory^m}\to \Sets^{\cat{M}}\]
is regular and conservative, but fails to preserve $\vee$ as well as $\bot$. Again we make explicit the least coverage $B$ forcing a standard interpretation. That is,  the least coverage such that the initial object $0$ is dense in $\evalmod{\syntob{}{P_{\bot}}}{-}$ and, for all disjunctions $\syntob{\alg{x}}{\phi\vee\psi}$ of \theory,  $\evalmod{\syntob{\alg{x}}{P_{\phi}}}{-}\vee\evalmod{\syntob{\alg{x}}{P_{\psi}}}{-}$ is dense in $\evalmod{\syntob{\alg{x}}{P_{\phi\vee\psi}}}{-}$. First, for all disjunctions(-in-context) $\syntob{\alg{x}}{\phi\vee\psi}$ of \theory, model diagrams $(D,F)$ in \cat{M}, and lists of elements $\alg{d}\in D^{l(\alg{x})}$, let $S_{\pair{\syntob{\alg{x}}{P_{\phi\vee\psi}},\alg{d}}}$ be the following cosieve on $(D,F)$:
\begin{align*}S_{\pair{\syntob{\alg{x}}{P_{\phi\vee\psi}},\alg{d}}}=&\cterm{\funktor{h}{(D,F)}{(D',F')}}{(D',F')\vDash P_{\phi}[ h(\alg{d})/\alg{x}]\vee P_{\psi}[ h(\alg{d})/\alg{x}]}
\end{align*}
%
%\begin{align*}S_{\pair{\syntob{\alg{x}}{P_{\phi\vee\psi}},\alg{d}}}=&\cterm{\funktor{h}{(D,F)}{(D',F')}}{(D',F')\vDash P_{\phi}[ h(\alg{d})/\alg{x}]}\\
%\cup & \cterm{\funktor{h}{(D,F)}{(D',F')}}{(D',F')\vDash P_{\psi}[ h(\alg{d})/\alg{x}]} 
%\end{align*}
%
Then let $B$ be specified by %adding to $E(D,F)$  the set of cosieves  
\[B(D,F)=E(D,F)\cup\cterm{S_{\pair{\syntob{\alg{x}}{P_{\phi\vee\psi}}},\alg{d}}}{\pair{\syntob{\alg{x}}{P_{\phi\vee\psi}},\alg{d}}\in F}\]
%
%where
%
%\begin{align*}S_{\phi\vee\psi,\alg{d}}=&\cterm{\funktor{h}{(D,F)}{(D',F')}}{h(\alg{d})\in \csem{\alg{x}}{P_{\phi}}^{(D',F')}}\\
%&\cup\cterm{\funktor{h}{(D,F)}{(D',F')}}{h(\alg{d})\in \csem{\alg{x}}{P_{\psi}}^{(D',F')}} \end{align*}
%
%(Here $h(\alg{d})\in\csem{\alg{x}}{P_{\phi}}^{(D',F')}$ is shorthand for\\
% $\cterm{\alg{c}\in D'}{h(\alg{d},\alg{c})}\subseteq\cterm{\alg{c}\in D'}{(D',F')\vDash \phi(\alg{c})}$.) 
%
%
Again,  $B$ is a coverage. We refer to it  as the \emph{minimal coverage}. A connection to Beth semantics will be displayed in Section \ref{subsubsection: beth completeness}. The proof of the following is similar to that of Theorem \ref{proposition: modified joyals theorem for coherent with C coverage}, and a corollary of it if \cat{M} is chase-complete, and is therefore omitted.

\begin{proposition}\label{proposition: modified joyals theorem for coherent}
Let \theory\ be an \textup{(}at least\textup{)} coherent theory, $\theory^m$ its regular Morleyization, \cat{M} a strongly conservative, small, full subcategory of $\moddiag{\theory^m}$, and $B$ be the minimal coverage on \op{\cat{M}}. Then the evaluation functor $\mo{Ev}$ composed with the sheafification functor $\funktor{a}{\Sets^{\cat{M}}}{\Sh{\op{\cat{M}},B}}$
\[\bfig
\qtriangle/>`>`@{<-^{)}}@<5pt>/<500,300>[\synt{C}{\theory}`\Sets^{\cat{M}}`\Sh{\op{\cat{M}},B};\operatorname{Ev}`a\circ\mo{Ev}`]
\qtriangle|abl|/>`>`{@{>}@<-5pt>}/<500,300>[\synt{C}{\theory}`\Sets^{\cat{M}}`\Sh{\op{\cat{M}},B};\operatorname{Ev}`a\circ\mo{Ev}`a]
\efig\]
is conservative, coherent,  and conditionally Heyting.
%
%Moreover, $\operatorname{Ev}:\synt{C}{\theory}\to<150>\Sh{\op{\cat{M}},E}$, in addition to being regular and conservative, $%%\operatorname{Ev}$ preserves the initial object and right adjoints to pullback functors whenever they exist.
%
\end{proposition}
%
%We refer to $B$ as the \emph{Beth coverage} and return to in Section \ref{section: veldman etc thms} where we %give a simpler description of it (and justify the choice of name for it).

If \cat{M} is chase-complete, then the minimal coverage on $\op{\cat{M}}$ can be strengthened while still yielding a conservative \theory-model as follows. 
Let $(D,F)\in \cat{M}$, let \syntob{\alg{x}}{\phi\vee \psi} be a disjunction of \theory, and let $\alg{d}\in D$ such that $\pair{\syntob{\alg{x}}{P_{\phi\vee\psi}},\alg{d}}\in F$. Then we have, 
\[\bfig
\Vtriangle/`<-`<-/<500,300>[\operatorname{Ch}(D,F\cup\{P_{\phi}{[}\alg{d}{]}\})`\operatorname{Ch}(D,F\cup\{P_{\psi}{[}\alg{d}{]}\})`(D,F);`c_0`c_1]
\efig\]
with $c_0$ the homomorphism induced by $(D,F)\subseteq (D,F\cup\{P_{\phi}{[}\alg{d}{]}\}) \subseteq \operatorname{Ch}(D,F\cup\{P_{\phi}{[}\alg{d}{]}\})$, and similarly for $c_1$.
We refer to such a pair, given by a fact of the form $\pair{\syntob{\alg{x}}{P_{\phi\vee\psi}},\alg{d}}$ in $F$, as a \emph{chase pair} 
%
%\pair{(D,F), \pair{\syntob{\alg{x}}{P_{\phi\vee\psi}},\alg{d}}}
%
over $(D,F)$. Since $\operatorname{Ch}$ is a functor, assigning to each $(D,F)\in\cat{M}$ the set of chase pairs over it is a coverage. To this we also add the coverage $E$, so that the family of covering families over $(D,F)$ is the union of the set of chase pairs and the set \cterm{\emptyset}{(D,F)\vDash P_\bot}. Denote the resulting coverage by $C$, and the least Grothendieck coverage containing $C$ by $\overline{C}$. Similarly, write  $\overline{B}$ for the least Grothendieck coverage containing $B$. We refer to $C$ as the \emph{disjunctive} coverage. 
%We display the following observations for reference.
%
The two coverages compare as follows.
%
%\begin{lemma}\label{lemma: B and C same if choice}
%Let \theory\ be a coherent theory, $\theory^m$ its regular Morleyization, and \cat{M} a chase-complete category of model diagrams for $\theory^m$. Let $\overline{B}$ and $\overline{C}$ be the least Grothendieck coverages containing the minimal and disjunctive coverages on \op{\cat{M}}, respectively. Then $\overline{B}\subseteq \overline{C}$. 
%If $\operatorname{Ch}$ is a weak reflection
%then also $\overline{C}\subseteq\overline{B}$.
%\end{lemma} 
%
\begin{lemma}\label{lemma: B and C same if choice}
Let \theory\ be a coherent theory, $\theory^m$ its regular Morleyization, and \cat{M} a chase-complete category of model diagrams for $\theory^m$. Let $\overline{B}$ and $\overline{C}$ be the least Grothendieck coverages containing the minimal and disjunctive coverages on \op{\cat{M}}, respectively. Then $\overline{B}\subseteq \overline{C}$. Moreover, 
\begin{enumerate}
\item If $\operatorname{Ch}$ is a weak reflection
then also $\overline{C}\subseteq\overline{B}$.
\item The statement that $\overline{C}= \overline{B}$ for arbitrary \theory\ is equivalent to the Axiom of Choice.
\end{enumerate} 

\begin{proof}

%
%For this,
%we refer to the proof of Theorem 5.16 in \cite{forsselllumsdaine} to the effect that if the regular theory with a single axiom $B(x)\vdash_x \fins{y}R(x,y)$ has a chase functor which is a weak reflection, then every surjection of sets splits. 

(1) and that $\overline{C}\supseteq \overline{B}$ is clear. Since (2) can be considered  more of  a remark  that will play no further role for us here, we only outline the proof:  
%We outline the proof that $\overline{C}= \overline{B}$ for arbitrary  \theory\ implies the Axiom of Choice.

Consider the \emph{coherent} theory \theory\ with  unary predicate symbols $B$ and $\bar{B}$,  one binary relation symbol $R$, and the axioms 
\begin{align*}
B(x)\wedge \bar{B}(x)\vdash_x \bot \quad \quad
%\top\vdash_x B(x)\vee \bar{B}(x) \quad \quad
B(x)\vdash_x \fins{y}R(x,y)
\end{align*}
A surjection $e:Y\to<125> X$ can be considered as a \theory-model, and therefore a $\theory^m$-model \alg{E} by, briefly, letting $|\alg{E}|$ be $X+Y$,  $X$ be the extension of $B$,  and the extension of  $R$ be the inverse of the graph of  $e$. 

Write $\beta^{\alg{x}}:= \bigwedge_{x\in\alg{x}}(B(x)\vee\bar{B}(x))$.
Let $(D,F)$ be the $\Sigma^m$-diagram presented by $D=\{\ast\}$ and\\  $F=\cterm{\pair{\syntob{\alg{x}}{P_{\phi}},\vec{\ast}}}{\syntob{\alg{x}}{\phi} \textnormal{ coherent in canonical context, and } \beta^{\alg{x}}\vdash^{\theory}_{\alg{x}}\phi}$.\\  
Then $(D,F)$ is a (finite) $\theory^m$-model diagram.  

Now, if $\overline{C}= \overline{B}$ then 
% by considering a basis for $\bar{C}$, that 
the family $S_{\pair{\syntob{x}{P_{B(x)\vee\bar{B}(x)}},{\ast}}}$ is the sieve generated by the chase pair given by \pair{\syntob{x}{P_{B(x)\vee\bar{B}(x)}},\ast}. Then, for a surjection $e:Y\to<125> X$ ,  elements in $X$ induces homomorphisms  $(D,F)\to<125>\alg{E}$, and the lifting through $\operatorname{Ch}(D,F\cup\{P_{B(x)}{[}\ast{]}\})$ gives a splitting of $e$ (see the proof of Theorem 5.16 in \cite{forsselllumsdaine}).

\end{proof}
\end{lemma} 
%
%
%ALTERNATIV: If the Axiom of Choice does not hold, then in general, it is not the case that $\overline{C}\subseteq \overline{B}$. 
%In spite of this, we have:
%

\begin{theorem}\label{proposition: modified joyals theorem for coherent with C coverage}\label{theorem: modified joyals theorem for coherent with C coverage}
Let \theory\ be an at least coherent theory, $\theory^m$ its regular Morleyization, \cat{M} a full, chase-complete subcategory of $\moddiag{\theory^m}$, and $C$ be the disjunctive coverage on \op{\cat{M}}. Then the evaluation functor $\mo{Ev}$ composed with the sheafification functor $\funktor{a}{\Sets^{\cat{M}}}{\Sh{\op{\cat{M}},C}}$
\[\bfig
\qtriangle/>`>`@{<-^{)}}@<5pt>/<500,300>[\synt{C}{\theory}`\Sets^{\cat{M}}`\Sh{\op{\cat{M}},C};\operatorname{Ev}`a\circ\mo{Ev}`]
\qtriangle|abl|/>`>`{@{>}@<-5pt>}/<500,300>[\synt{C}{\theory}`\Sets^{\cat{M}}`\Sh{\op{\cat{M}},C};\operatorname{Ev}`a\circ\mo{Ev}`a]
\efig\]
is conservative, coherent,  and conditionally Heyting.
\begin{proof} Let $U^n=\evalmod{\syntob{\alg{x}}{\top}}{-}$ where the length of \alg{x} is $n$. Recall from e.g.\ \cite{maclane92} that the (coherent) sheafification functor $a$ restricts to (one half of) a poset isomorphism between  closed subobjects of $U^n$ and subobjects of $a(U^n$),
\[\operatorname{ClSub}_{}(U^n)\cong\operatorname{Sub}_{}(a(U^n))\]
It follows that it is sufficient to  show that, for any formula $\syntob{\alg{x}}{\theta}$ of \theory, the functor $\evalmod{\syntob{\alg{x}}{P_{\theta}}}{-}$ is $C$-closed in the subobject lattice of $U^{n}$. Let  $(D,F)\in \cat{M}$, let \syntob{\alg{y}}{\phi\vee\psi} be a disjunction of \theory, and let $\alg{d}\in D^{\operatorname{length}(\alg{y})}$ such that $(D,F)\vDash P_{\phi\vee\psi}[\alg{d}/\alg{y}]$. Let $(D_1,F_1)= \operatorname{Ch}(D,F\cup\{P_{\phi}[\alg{d}]\}))$ and $(D_2,F_2)= \operatorname{Ch}(D,F\cup\{P_{\psi}[\alg{d}]\}))$. Let $\alg{c}\in D^n$ and assume $(D_1,F_1),(D_2,F_2)\vDash P_{\theta}[\alg{c}/\alg{x}]$. Then $\theory^m_{(D,F)}$ proves the sequents $P_{\phi}[\alg{d}/\alg{y}]\vdash P_{\theta}[\alg{c}/\alg{x}]$ and $P_{\psi}[\alg{d}/\alg{y}]\vdash P_{\theta}[\alg{c}/\alg{x}]$. Thus there exists  proofs with premisses in $\theory^m$ of  $\chi\wedge P_{\phi}[\alg{d}/\alg{y}]\vdash P_{\theta}[\alg{c}/\alg{x}]$ and $\xi\wedge P_{\psi}[\alg{d}/\alg{y}]\vdash P_{\theta}[\alg{c}/\alg{x}]$ where $\chi$ and $\xi$ are conjunctions of atomic  sentences over $\Sigma\cup D$ which are true in $(D,F)$. By Lemma \ref{lemma: replacing constants with variables}, there are proofs with the same premisses of  sequents $\chi'\wedge P'_{\phi}\vdash_{\alg{w}} P'_{\theta}$ and $\xi'\wedge P''_{\psi}\vdash_{\alg{v}} P''_{\theta}$ over $\Sigma$, where we assume $\alg{w}$ and $\alg{v}$ disjoint, and a function \funktor{f}{\alg{w},\alg{v}}{D} such that $\chi'[f]=\chi$ and $\xi'[f]=\xi$ and $P'_{\phi}[f]=P_{\phi}[\alg{d}/\alg{y}]$ and $P''_{\psi}[f]=P_{\psi}[\alg{d}/\alg{y}]$ and $P'_{\theta}[f]=P''_{\theta}[f]=P_{\theta}[\alg{c}/\alg{x}]$. 
Form the finite set $E$ of equalities with variables from the (assumed disjoint) lists  \alg{w},\alg{v}, \alg{x}, and \alg{y} as follows: for each argument slot in $\alg{P}_{\phi}$, with that slot occupied by, say, $y$ in $P_{\phi}$, $w$ in $P'_{\phi}$, and $v$ in $P''_{\phi}$, add $y=w$ and $y=v$ to $E$; for each argument slot in $\alg{P}_{\psi}$, with that slot occupied by, say, $y$ in $P_{\psi}$, $w$ in $P'_{\psi}$, and $v$ in $P''_{\psi}$, add $y=w$ and $y=v$ to $E$; and for each argument slot in $\alg{P}_{\theta}$, with that slot occupied by, say, $x$ in $P_{\theta}$, $w$ in $P'_{\theta}$, and $v$ in $P''_{\theta}$, add $x=w$ and $x=v$ to $E$. Let $\rho$ be the conjunction of the equalities in $E$. Extend $f$ by: if $x$ occurs in an equality in $\rho$, say $x=w$, then add $x\mapsto f(w)$; and if $y$ occurs in an equality in $\rho$, say $y=w$, then add $y\mapsto f(w)$. Notice that this is well defined and $P_{\phi}[f]=P_{\phi}[\alg{d}/\alg{y}]$, $P_{\psi}[f]=P_{\psi}[\alg{d}/\alg{y}]$, and $P_{\theta}[f]=P_{\theta}[\alg{c}/\alg{x}]$. Then: 1) $\rho[f]$ is true in $(D,F)$ and
 2) the sequents
$\fins{\alg{w},\alg{v}}\rho \wedge \chi' \wedge \xi' \wedge P_{\phi}\vdash_{\alg{x},\alg{y}}P_{\theta}$
and
$\fins{\alg{w},\alg{v}}\rho \wedge \chi' \wedge \xi' \wedge P_{\psi}\vdash_{\alg{x},\alg{y}}P_{\theta}$
are provable in $\theory^m$. We can replace the regular formula $\fins{\alg{w},\alg{v}}\rho \wedge \chi' \wedge \xi' $ by an atomic formula in $\Sigma^m$; say $P_{\gamma}$. We then have that \theory\ proves the sequents   
$\gamma \wedge {\phi}\vdash_{\alg{x},\alg{y}}{\theta}$
and
$\gamma \wedge {\psi}\vdash_{\alg{x},\alg{y}}{\theta}$. Thus \theory\ proves $\gamma \wedge (\phi\vee \psi)\vdash_{\alg{x},\alg{y}}{\theta}$, whence $\theory^m$ proves $P_{\gamma} \wedge P_{\phi\vee \psi}\vdash_{\alg{x},\alg{y}}P_{\theta}$. Since substituting \alg{d} for \alg{y} and \alg{c} for \alg{x} makes the antecedent true in the $\theory^m$-model $(D,F)$, $P_{\theta}[\alg{c}/\alg{x}]$ must also be true in $(D,F)$.  
%
%
%With $\evalmod{\syntob{\alg{x}}{P_{\theta}}}{-}$  $C$-closed in the subobject lattice of $U^n$, and $\operatorname{Ev}$ conditionally sub-Heyting, we then have the following with respect to the Heyting operations in the lattices of closed subobjects of the presheaves $U^n$. First, as in Lemma \ref{lemma: modified joyals theorem for regular},  $\evalmod{\syntob{\alg{x}}{P_{\bot}}}{-}$ is the initial closed subobject, ${0}_{cl}$. For the meet, $\evalmod{\syntob{\alg{x}}{P_{\theta}}}{-}\wedge_{cl} \evalmod{\syntob{\alg{x}}{P_{\phi}}}{-}=\evalmod{\syntob{\alg{x}}{P_{\theta}}}{-}\wedge \evalmod{\syntob{\alg{x}}{P_{\phi}}}{-}=\evalmod{\syntob{\alg{x}}{P_{\theta\wedge\phi}}}{-}$. For the join we have $\evalmod{\syntob{\alg{x}}{P_{\theta}}}{-}\vee_{cl} \evalmod{\syntob{\alg{x}}{P_{\phi}}}{-}=\evalmod{\syntob{\alg{x}}{P_{\theta\vee\phi}}}{-}$ since $C$ extends $B$. For the image, we have $\exists_{cl} \evalmod{\syntob{\alg{x,\alg{y}}}{P_{\psi}}}{-}=\evalmod{\syntob{\alg{y}}{P_{\fins{\alg{x}}\psi}}}{-}$ since $\operatorname{Ev}$ is regular and  $\evalmod{\syntob{\alg{y}}{P_{\fins{\alg{x}}\psi}}}{-}$ is closed. And, by the conditional sub-Heytingness of $\operatorname{Ev}$, we have  $\evalmod{\syntob{\alg{x}}{P_{\theta}}}{-}\rightarrow_{cl} \evalmod{\syntob{\alg{x}}{P_{\phi}}}{-}=\evalmod{\syntob{\alg{x}}{P_{\theta}}}{-}\rightarrow \evalmod{\syntob{\alg{x}}{P_{\phi}}}{-}=\evalmod{\syntob{\alg{x}}{P_{\theta\rightarrow\phi}}}{-}$ and    $\forall_{cl} \evalmod{\syntob{\alg{x,\alg{y}}}{P_{\psi}}}{-}=\evalmod{\syntob{\alg{y}}{P_{\alle{\alg{x}}\psi}}}{-}$.   
%
%
\end{proof}
\end{theorem}
%
%Proposition \ref{proposition: modified joyals theorem for coherent} now follows from the proof of Proposition \ref{proposition: modified joyals theorem for coherent with C coverage} by noting that the evaluation functors must be $B$-closed since they are $C$-closed and $C$ extends $B$. 
% 

\subsubsection{Kripke and generalized Beth models}\label{subsection: models on posets}%%%%%%%%%%%%%%%%%%%%%%%%%%%%%%%%%%%%%%%%%%%%%%%%%%%%%%%%%%%%%%%%%%%

The models of  \ref{proposition: modified joyals theorem for regular} and  \ref{proposition: modified joyals theorem for coherent with C coverage} can be translated to models in presheaves and sheaves on posets using e.g.\ the Diaconescu cover (see e.g.\ \cite{maclane92} for a description of the Diaconescu cover). However, in our current setting, we can use, more directly, the poset of model diagrams and inclusions. 
 We state this also as covering lemma, the technical heart of which is the following.  Write  $\operatorname{MDiag}^{\subseteq}(\theory)$ for the poset of of model diagrams and inclusions.  Write  $\funktor{\pi}{\operatorname{MDiag}^{\subseteq}(\theory)}{\moddiag{\theory}}$ for the functor that sends an inclusion to the homomorphism it induces.

%The inclusion induces, when restricted to suitable sets of models, an open surjection of respective presheaf toposes, in virtue of satisfying the following technical condition (\emph{cf}.\  \cite[C3.1.2]{elephant1}).

%The following lemma verifies the condition of \cite[?]{elephant1}

\begin{lemma}\label{lemma: open cover}
Let  \moddiag{\theory} be the category of model diagrams for some (regular or coherent) theory \theory. 
For any homomorphism $\funktor{h}{(D_0,F_0)}{(D_1,F_1)}$ where $D_0$ and $D_1$ are disjoint there exists an extension $a\!:\!(D_0,F_0)\subseteq (D_2,F_2)$ and homomorphisms $r$ and $i$ as in the following diagram
\[\bfig
\dtriangle|arb|/<-`@{>}@/^5pt/`>/<750,500>[(D_1,F_1)`(D_0,F_0)`(D_2,F_2);h`i`\pi(a)]
\dtriangle|alb|/<-`@{<-}@/_5pt/`>/<750,500>[(D_1,F_1)`(D_0,F_0)`(D_2,F_2);`r`\pi(a)]
%\place(375,0)[\subseteq]
\efig\]
such that the outer triangle commutes and $r\circ i$ is the identity on $(D_1,F_1)$.
\begin{proof} Let  $D_2=D_0 \cup D_1$ and let $F_2$ be the diagram generated by $F_0\cup F_1\cup \cterm{d=d'}{h(d,d')}$. We have
 $a:(D_0,F_0)\subseteq (D_2,F_2)$. 
%
%Explicitly, $a(d,d')\Leftrightarrow (d=d')\in F_2$. Similarly, 
%
Let  $i$ be the homomorphism induced by the inclusion $(D_1,F_1)\subseteq (D_2,F_2)$. Then
\[\pi(a)(d,d') \Leftrightarrow \fins{c\in D_1} h(d,c)\wedge i(c,d')\]
so the outer triangle commutes. Let $r(d',d)\Leftrightarrow (d'=d)\in F_1 \vee h(d',d)$. It is then  straightforward that $r$ is a well-defined homomorphism, as well as both a left and right inverse to $i$, from which it also follows that $(D_2,F_2)$ is a \theory-model diagram.
\end{proof}  
\end{lemma}
Clearly, the assumption that the domain and codomain of $h$ are disjoint can be done without loss if the codomain can be replaced by an isomorphic copy disjoint from both it and the domain. 
If \cat{M} is a chase-complete set of model diagrams for a regular theory \theory\ and $\cat{M}^{\subseteq}$ is the poset of model diagrams in \cat{M} and inclusions, then Lemma \ref{lemma: open cover} implies (by e.g.\ \cite[C3.1.2]{elephant1}) that the restriction functor 
\[\Sets^{\cat{M}_{}}\to^{\pi^*}\Sets^{\cat{M}_{}^{\subseteq}}\]
is Heyting and conservative. By Theorem \ref{theorem: joyals theorem} we then obtain the following new version of Joyal's theorem:
\begin{theorem}\label{Theorem: joyals with inclusions}
Let \theory\ be a regular theory and $\cat{M}_{}^{\subseteq}$ a chase-complete set of model diagrams ordered by inclusions. Then the evaluation functor 
\[\synt{C}{\theory}\to \Sets^{\cat{M}^{\subseteq}_{}}\]
is conservative and conditionally sub-Heyting.
\end{theorem}
%

%As we began this section by saying, 
We have, as corollaries to  Theorem \ref{theorem: joyals theorem} (or Proposition \ref{proposition: modified joyals theorem for regular}), (fallible) Kripke completeness results for  theories in certain fragments of first-order logic.
%; where ``general'' refers to the absence of  assumptions regarding size, dicreteness, or the derivability of the sequent $(\top\vdash\fins{x}\top)$.
%
By \ref{lemma: open cover},
the underlying poset of the Kripke models can be taken to be a set of model diagrams for the regular Morleyization of \theory\ ordered by inclusion. 
Similarly, as a corollary of Theorem \ref{proposition: modified joyals theorem for coherent with C coverage} we have a completeness theorem for first-order theories with respect to a generalized version of Beth semantics. 
These are  fairly straighforward cases of translating models on presheaves and sheaves on posets to Kripke and Beth-style presentations. For explicitness, we give some further details, also making clear  
%We write these corollaries out here for explicitness, also making explicit 
what notions of Kripke and ``generalized'' Beth models we have in mind:

%For explicitness, we display  these corollaries with respect to the following notion of Kripke and generalized Beth model.

Let $\Sigma$ be a relational signature. Write $\topo{S}$ for  the partially ordered class of $\Sigma$-diagrams and (homomorphic) inclusions. Thus an object $S$ in $\topo{S}$ is a Tarski structure for $\Sigma$ with a congruence relation $\sem{=}^S$ interpreting $=$.    Write $\topo{F}$ for  the partially ordered class of \emph{fallible} $\Sigma$-diagrams and diagram inclusions:  an object $S$ in $\topo{F}$ is an inhabited diagram for $\Sigma$ with a subset $\sem{\bot}^S\subseteq 1$ of the terminal set interpreting $\bot$, and such that $S$ satisfies the axioms $\bot \vdash_{\alg{x}} \phi$ for all atomic formulas $\phi$ in canonical context \alg{x} over $\Sigma$. The inclusions in \topo{F} must preserve $\sem{\bot}$, i.e.\ $S\subseteq S' \Rightarrow \sem{\bot}^S\subseteq \sem{\bot}^{S'}$.

%By ``$\Sigma$-structure'' in the following definition we mean what we elsewhere have been referring to as a diagram, i.e.\ an ordinary (but possibly empty) Tarski structure  in \Sets\ for $\Sigma$ where equality is interpreted as a congruence. By an ``fallible $\Sigma$-structure'' we mean a $\Sigma$-structure where, additionally, $\bot$ is interpreted as an arbitrary truth value---i.e.\ $\bot$ is given an interpretation as if it were a propositional variable---and the structure satisfies the axioms $\bot\vdash_{\alg{x}}\phi$, for all atomic formulas-in-context \FIC{\alg{x}}{\phi}, \emph{and} the axiom $\bot\vdash \fins{x}\top$. A homomorphism between fallible structures must preserve the interpretation of $\bot$ (as if it were a propositional variable).
%
\begin{definition}\label{definition: generalized beth model}

(I)
Let $\Sigma$ be a relational signature. 
By a \emph{generalized \textup{(}fallible\textup{)} Beth structure} for $\Sigma$ we mean a triple \pair{P,D,T} where $P$ is a poset; $\mathfrak{D}$ is a functor from  $P$ to \topo{S} (\topo{F}); and $T$ is an assignment of inhabited sets of subposets of $P$ to nodes of $P$ such that:
\begin{enumerate}[(i)]
\item all elements of $T(p)$ are finite, binary trees with root $p$, and $T(p)$ is closed under initial binary subtrees;
\item if $t\in T(p)$ with leaf nodes $q_1,\ldots,q_n$ and $t_1\in T(q_1),\ldots, q_n\in T(q_n)$ then the tree obtained by extending $t$ with the $t_i$'s   is in $T(p)$;
and for all $p\in P$, $t\in T(p)$, and $q\in t$, $t\cap \uparrow q \in T(q)$;
 and
\item for all  $p\leq p'$ in $P$ and $t\in T(p)$ there exists $t'\in T(p')$ such that for all leaf nodes $q'$ of $t'$ there exists a leaf node $q$ of $t$ such that $q\leq q'$.
\end{enumerate}
The clauses of the forcing relation $p\Vdash \phi[\alg{d}/\alg{x}]$ between $p\in P$, first-order formulas-in-context \FIC{\alg{x}}{\phi}, and $\alg{d}\in \mathfrak{D}(p)^{l(\alg{x})}$ are then: 
\begin{enumerate}[(a)]
\item for $\phi$ atomic or equal $\bot$ or $\top$, $p\Vdash \phi[\alg{d}/\alg{x}]$ if there exists $t\in T(p)$ such that for all leaf nodes $q\in t$ it is the case that $\mathfrak{D}(q)\vDash \phi[\alg{d}/\alg{x}]$;
\item for $\phi=\psi\wedge\theta$,  $p\Vdash \phi[\alg{d}/\alg{x}]$ if $p\Vdash \psi[\alg{d}/\alg{x}]$ and $p\Vdash \theta[\alg{d}/\alg{x}]$;
\item for $\phi=\psi\vee\theta$,  $p\Vdash \phi[\alg{d}/\alg{x}]$ if there exists $t\in T(p)$ such that for all leaf nodes $q$ in $t$ it is the case that $q\Vdash \psi[\alg{d}/\alg{x}]$ or $q\Vdash \theta[\alg{d}/\alg{x}]$;
\item for $\phi=\psi\rightarrow\theta$,  $p\Vdash \phi[\alg{d}/\alg{x}]$ if for all $p'\geq p$ it is the case that if $p'\Vdash \psi[\alg{d}/\alg{x}]$ then $p'\Vdash \theta[\alg{d}/\alg{x}]$;
\item for $\phi=\fins{y}\psi$,  $p\Vdash \phi[\alg{d}/\alg{x}]$ if  there exists $t\in T(p)$ such that for all leaf nodes $q$ in $t$ there exists $c\in \mathfrak{D}(q)$ such that $q\Vdash \psi[c/y,\alg{d}/\alg{x}]$; and 
\item for $\phi=\alle{y}\psi$,  $p\Vdash \phi[\alg{d}/\alg{x}]$ if for all $p'\geq p$ and all $c\in \mathfrak{D}(p')$ it is the case that $q\Vdash \psi[c/y,\alg{d}/\alg{x}]$.
\end{enumerate}

(II) By a \emph{ \textup{(}fallible\textup{)} Kripke structure} we mean  a generalized (fallible) Beth structure where $T(p)$ contains only the one node tree on $p$. 

(III) By  (\emph{fallible}) \emph{Beth structure} we mean a generalized (fallible) Beth structure where $P$ is a binary tree 
 and $T(p)$ is the set of initial binary subtrees of $\uparrow(p)$. (This notion of Beth structure is, then, with respect to the strong rather than the weak notion of forcing, cf.\ \cite[Ch.13 1.8]{troelstravandalen:88ii}.)

(IV) By a (\emph{generalized, fallible}) \emph{Beth$^{\star}$ structure} we mean a (generalized, fallible) Beth structure where ``covers are only relevant for disjunctions'', i.e.\ one satisfying the following additional conditions:
\begin{enumerate}[(1)]
\item for $\phi$ atomic or $\bot$, it is the case that $p\Vdash\phi[\alg{d}]$ iff $p\vDash\phi[\alg{d}]$, and
\item for all formulas of the form $\fins{x}\phi\in \cat{L}$, it is the case that  $p\Vdash\fins{x}\phi[\alg{d}]$ iff there exists $a\in p$ such that $p\Vdash\phi[a,\alg{d}]$.  
\end{enumerate} 
\end{definition}

We state corollaries of Theorem \ref{Theorem: joyals with inclusions}, Proposition \ref{proposition: modified joyals theorem for regular}, and Theorem \ref{proposition: modified joyals theorem for coherent with C coverage} in terms of Definition \ref{definition: generalized beth model}. Let the \emph{$\vee$-free} fragment of FOL be the fragment consisting of sequents not mentioning the connective $\vee$, and the \emph{$\bot$,$\vee$-free} fragment be the one not mentioning $\bot$ or $\vee$.
 \begin{corollary}\label{corollary: Kripke comp for bot and disj free}
Let \theory\ be a   theory in the $\bot$,$\vee$-free fragment over the signature $\Sigma$. Then there exists a Kripke model for \theory\ which is conservative (with respect to the $\bot$,$\vee$-free fragment).
\begin{proof}
Let $\theory^m$ over $\Sigma^m$ be the regular Morleyization of \theory. By Theorem \ref{theorem: Reg theories have replete sets of models} there exists a chase-complete category \cat{M} of model diagrams for  $\theory^m$. By Theorem \ref{Theorem: joyals with inclusions}, the evaluation functor  $\operatorname{Ev}:\synt{C}{\theory}\simeq\synt{C}{\theory^m}\to \Sets^{\cat{M}^{\subseteq}_{}}$ is conditonally sub-Heyting, thus giving a conservative (w.r.t.\ the $\bot$,$\vee$-free fragment) model of \theory. Define a Kripke structure $K$ by letting the poset $P$ be $\cat{M}^{\subseteq}$ and  the functor $\funktor{\mathfrak{D}}{\cat{M}^{\subseteq}}{\topo{S}}$ be the forgetfull functor. Let \FIC{\alg{x}}{\phi} be $\bot$,$\vee$-free over $\Sigma$, $S\in \cat{M}^{\subseteq}$, and $\alg{d}$ a list of the same length as \alg{x} of elements  in the domain of $S$ . Using that $S\vDash P_{\phi}[\alg{d}/\alg{x}]\Leftrightarrow [\alg{d}]\in \evalmod{\FIC{\alg{x}}{\phi}}{S}$ and that $\operatorname{Ev}$ is conditionally sub-Heyting, a straightforward induction argument on $\FIC{\alg{x}}{\phi}$ shows that $S\Vdash^{K} \phi[\alg{d}/\alg{x}]\Leftrightarrow S\vDash P_{\phi}[\alg{d}/\alg{x}] $.   

\end{proof}
\end{corollary}
 \begin{corollary}\label{corollary: Mod Kripke comp for disj free}
Let \theory\ be a $\vee$-free theory. Then there exists a fallible Kripke model for \theory\ which is conservative with respect to the \emph{$\vee$-free} fragment.
\begin{proof}
From Proposition \ref{proposition: modified joyals theorem for regular} and Theorem \ref{Theorem: joyals with inclusions} (similarly to  \ref{corollary: Kripke comp for bot and disj free}).
% cf.\ also the proof of Corollary \ref{corollary: generalized beth model}). 

%From Proposition \ref{proposition: modified joyals theorem for regular} (similarly to how \ref{corollary: Kripke comp for bot and disj free} is obtained from \ref{proposition: modified joyals theorem for regular}). 
%
%\[\bfig  \efig\]
%
\end{proof}
\end{corollary}
\begin{remark}\label{Remark: completeness for disjunction free}
%Note that the restriction to the  $\bot$,$\vee$-fragment is essential in  \ref{corollary: Kripke comp for bot and disj free}. 
%
The restrictions are essential.
 The existence of a conservative Kripke model (that is, a non-fallible one) for $\vee$-free theories implies LEM (see \cite{mccarty:08})\footnote{In fact, it is equivalent to it, since with LEM we can  distinguish between exploding and non-exploding models.}.
%; and LEM implies the existence of such a model by Propostion \ref{proposition: modified joyals theorem for regular} and  the observation that LEM implies the equivalence
%\[\Sh{\op{\cat{M}}, \thry{E}}\simeq \Sets^{\cat{M}^{st}}\]
%where $\cat{M}^{st}\to/^{(}->/<125> \cat{M}$ is the subcategory of standard  (non-exploding) models.
%
%On the other hand, 
A Kripke completeness theorem for $\bot$-free theories, or a fallible Kripke completeness theorem for full FOL, would imply e.g.\ that the Boolean Prime Ideal theorem is provable in ZF.
For theories whose axioms do not mention $\bot$ or $\vee$, such as  the  empty theory, the existence of a Kripke model which is conservative with respect to all first-order sequents implies MP (see \cite{mccarty:08}).      
%
%Note, further, that a completeness theorem for modified Kripke semantics for full first order logic cannot be proven in IZF.
%Since, classically, modified Kripke semantics is just the usual Kripke semantics, the completeness of the former is classically equivalent to the completeness of the latter, and Kripke completeness for first-order theories is, of course, not provable in ZF.

% Hence, if MKC was derivable in IZF, we would have ZF=IZF+LEM=IZF+MKC+LEM=ZF+KC. which is absurd since ZF+KC implies BPI, %which is known to be independent of ZF (see, e.g., \cite{jech}).
\end{remark}
As an example application of Theorem \ref{Theorem: joyals with inclusions}, we give a short,  semantic proof of the disjunction property for (arbitrary) $\vee$-free theories (cf.\ \cite{troelstraschwichtenberg:96}) by reducing it to the disjunction property for regular theories (see e.g.\ \cite{elephant1}; note that the disjunction property for regular theories  also directly follows from Proposition \ref{theorem: chase completeness for regular theories}).
%; if $\phi$, $\psi$, and $\theta$ are regular formulas and $\phi\vdash_{\alg{x}}\psi\vee \theta$ is provable from regular axioms, then so is either $\phi\vdash_{\alg{x}}\psi$ or $\phi\vdash_{\alg{x}} \theta$.)

\begin{corollary}\label{corollary: disjunction property}
Let \theory\ be a first-order theory the axioms of which are $\vee$-free, and let $\phi$, $\psi$, and $\theta$ be $\vee$-free formulas. If \theory\ proves the sequent   $\phi\vdash_{\alg{x}}\psi\vee \theta$, then \theory\ proves either $\phi\vdash_{\alg{x}}\psi$ or $\phi\vdash_{\alg{x}} \theta$.
\begin{proof}
Consider the fallible Kripke model 
%$\synt{C}{\theory}\simeq \synt{C}{\theory^m}\to \Sets^{\cat{M}^{\subseteq}_{}}$ of   Theorem \ref{Theorem: joyals with inclusions}. 
of Corollary \ref{corollary: Mod Kripke comp for disj free}. (As in \ref{corollary: Kripke comp for bot and disj free}) the nodes are  models of $\theory^m$ and for all $\vee$-free formulas $\xi$ we have $(D,F)\vDash P_{\xi}[\alg{d}/\alg{x}]$ $\Leftrightarrow$ $(D,F)\Vdash \xi[\alg{d}/\alg{x}]$. Then 
for any $(D,F)$ in $\cat{M}^{\subseteq}$ and $\alg{d}\in D^{l(\alg{x})}$ we have: $(D,F)\vDash P_{\phi}[\alg{d}/\alg{x}]$ $\Leftrightarrow$ $(D,F)\Vdash \phi[\alg{d}/\alg{x}]$ $\Rightarrow$ $(D,F)\Vdash (\psi\vee \theta)[\alg{d}/\alg{x}]$  $\Rightarrow$ $(D,F)\Vdash \psi[\alg{d}/\alg{x}]$ or $(D,F)\Vdash \theta[\alg{d}/\alg{x}]$ $\Leftrightarrow$ $(D,F)\vDash P_{\psi}[\alg{d}/\alg{x}]$ or $(D,F)\vDash P_{\theta}[\alg{d}/\alg{x}]$ $\Leftrightarrow$ $(D,F)\vDash (P_{\psi}\vee P_{\theta})[\alg{d}/\alg{x}]$. Thus by Lemma \ref{lemma: chasecomplete gives geometric conservative},  $\theory^m$ proves the sequent 
%
% and $P_{\varphi}(\alg{x})$ correspondingly in $\Sigma^m$, that $(D,F)\Vdash\varphi[\alg{d}/\alg{x}] \Leftrightarrow (D,F)\vDash P_{\varphi}[\alg{d}/\alg{x}]$. 
%
%But, noting the forcing condition for disjunction in Kripke models, this means that the coherent sequent 
%
$P_{\phi}\vdash_{\alg{x}} P_{\psi}\vee P_{\theta}$. Therefore, $\theory^m$ proves the sequent $P_{\phi}\vdash_{\alg{x}} P_{\psi}$ or the sequent $P_{\phi}\vdash_{\alg{x}} P_{\theta}$, whence \theory\ proves  $\phi\vdash_{\alg{x}}\psi$ or $\phi\vdash_{\alg{x}} \theta$.

% is valid with respect to the $\theory^{m}$-models in \cat{M},  and therefore provable in the regular theory $\theory^m$.  

%$(D,F)\in \cat{M}^{\subseteq}$, 
%
%Let $(D,F)$ and \alg{d} be given, and assume $(D,F)\Vdash\phi[\alg{d}/\alg{x}]$. Then $(D,F)\Vdash(\psi\vee\theta}[\alg{d}/\alg{x}]$, and thus either $(D,F)\Vdash\psi[\alg{d}/\alg{x}]$ or $(D,F)\Vdash\theta[\alg{d}/\alg{x}]$. But this means that the sequent $P_{\phi}(\alg{x})\vdash_{\alg{x}} P_{\psi}(\alg{x})\vee P_{\theta}$   
\end{proof}
\end{corollary}
Finally, for full first-order logic we have:
 
\begin{corollary}\label{corollary: generalized beth model}
Let \theory\ be a first-order theory. Then \theory\ has a conservative generalized fallible Beth$^{\star}$ model.
\begin{proof}
Let $\theory^m$ be the regular Morleyization of \theory\ over extended signature $\Sigma^m$, with $\Sigma$ the signature of \theory. By Theorem \ref{theorem: Reg theories have replete sets of models} there exists a chase-complete category \cat{M} of model diagrams for  $\theory^m$. By Theorem \ref{proposition: modified joyals theorem for coherent with C coverage} the functor $a\circ \operatorname{Ev}:\synt{C}{\theory}\simeq\synt{C}{\theory^m}\to \Sh{\op{\cat{M}},  C}$ is conservative and Heyting. From Lemma \ref{lemma: open cover}, by \cite[C2.3.18--19(i)]{elephant1} and \cite[C3.1.23]{elephant1},  the (right) top functor  of the following commutative (up to isomorphism) diagram
\[\bfig 
\square(1000,0)/>`<-`<-`>/<1000,500>[\Sh{\op{\cat{M}},  C}`\Sh{\op{\cat{M}^{\subseteq}},C}`\Sets^{\cat{M}}`\Sets^{\cat{M}^{\subseteq}};a\circ\pi^{\ast}\circ i`a`a`\pi^{\ast}]

\dtriangle/<-`<-`>/<1000,500>[\Sh{\op{\cat{M}},  C}`\synt{C}{\theory}\simeq\synt{C}{\theory^m}`\Sets^{\cat{M}};a\circ \operatorname{Ev}``\operatorname{Ev}]
  \efig\]
is Heyting and conservative. Hence so is the composite top  functor. As in the proof of \ref{proposition: modified joyals theorem for coherent with C coverage}, the subpresheaves of the form $\pi^{\ast}\circ\operatorname{Ev}_{\syntob{\alg{x}}{\phi}}\to/^{(}->/<125> \pi^{\ast}\circ\operatorname{Ev}_{\syntob{\alg{x}}{\top}}$ are $C$-closed.

Define a generalized fallible Beth structure $B$ as follows. Let $P$ be $\cat{M}^{\subseteq}$, and let the functor $\funktor{\mathfrak{D}}{\cat{M}^{\subseteq}}{\topo{F}}$ send a $\Sigma^m$-diagram $S$ to its $\Sigma$ reduct extended with $\sem{\bot}^S:=\sem{P_{\bot}}^S$.
For $S\in\cat{M}^{\subseteq}$ let $T(S)$ be the set of finite binary trees with nodes in $\cat{M}^{\subseteq}$,  root $S$, and such that the children of any node $S'$ form a chase pair  (as given in the paragraph following \ref{proposition: modified joyals theorem for coherent}) over $S'$. 

Let \FIC{\alg{x}}{\phi} be first-order over $\Sigma$, $S\in \cat{M}^{\subseteq}$, and $\alg{d}$ a list of the same length as \alg{x} of elements  in the domain of $S$ .
We show by induction on $\phi$ that $S\Vdash^{B} \phi[\alg{d}/\alg{x}]\Leftrightarrow S\vDash P_{\phi}[\alg{d}/\alg{x}] $. 
And, simultaneously, that $B$ satisfies the conditions for being a Beth$^{\star}$ structure.

Let $\phi$ be atomic or $\phi = \bot$. Suppose $S\Vdash^{B} \phi[\alg{d}/\alg{x}]$. Then there exists a tree $t\in T(S)$ such that for all leaves $S'$ it is the case that $\mathfrak{D}(S')\vDash \phi[\alg{d}]$. Hence  $S'\vDash P_{\phi}[\alg{d}/\alg{x}] $. Now, the inclusions $S\subseteq S'$ define a $C$-cover, so  $S\vDash P_{\phi}[\alg{d}/\alg{x}] $, and thus $\mathfrak{D}(S)\vDash \phi[\alg{d}/\alg{x}] $. The converse is immediate. The case for the existential quantifier is similar, and the the case for conjuction is immediate. 

Let  $\phi = \psi \vee \theta $.   Suppose $S\Vdash^{B} \phi[\alg{d}/\alg{x}]$. Then there exists a tree $t\in T(S)$ such that for all leaves $S'$ it is the case that $S'\Vdash^{B}  \psi[\alg{d}]$ or $S'\Vdash^{B}  \psi[\alg{d}]$. By induction hypothesis  $S'\vDash P_{\psi}[\alg{d}]$ or $S'\vDash P_{\theta}[\alg{d}]$, so $S'\vDash P_{\phi}[\alg{d}]$, and since $t$ defines a $C$-cover on $S$, $S\vDash P_{\phi}[\alg{d}]$. Conversely, suppose $S\vDash P_{\phi}[\alg{d}/\alg{x}]$. Then \pair{\syntob{\alg{x}}{\phi},\alg{d}} defines a chase pair on $S$, yielding a tree in $T(S)$ with two leaves $S'$ and $S''$ such that, by the induction hypothesis, $S'\Vdash^B {\psi}[\alg{d}]$ and $S''\Vdash^B {\phi}[\alg{d}]$.

Let  $\phi = \alle{y}\psi$.   Suppose $S\Vdash^{B} \phi[\alg{d}/\alg{x}]$.   Then for all $S\subseteq S'$ we have   $\mathfrak{D}(S')\Vdash^{B} \psi[\alg{d}/\alg{x},c/y]$ for all $c$ in the domain of $S'$, thus by induction hypothesis $S'\vDash P_{\psi}[\alg{d},c]$. Hence, since $\funktor{\pi^{\ast}\circ \operatorname{Ev}}{\synt{C}{\theory}\simeq \synt{C}{\theory^m}}{\Sets^{\cat{M}^{\subseteq}}}$ is Heyting, $S\vDash P_{\phi}[\alg{d}/\alg{x}] $. The converse follows from $P_{\alle{y}\psi}\vdash_{\alg{x},y}^{\theory^m}P_{\psi}$. And $\rightarrow$ is similar.

\end{proof}
\end{corollary}

\subsubsection{Beth completeness}\label{subsubsection: beth completeness}
%
%

%We regard the site $(\op{\cat{M}},C)$ 
%(and to some extent also $(\cat{M},B)$)
%as a ``Beth-like'' or a ``generalized Beth'' model. It consists of a category (or poset, as we return to in  Section \ref{subsection: models on posets}) where every object/node comes equipped with a family of finite binary trees with respect to which the forcing condition for disjunction is specified. Moreover, in the enumerable\ case we can derive a version of (fallible) Beth completeness from it. In the following we take the notion of Beth structure to be with respect to the \emph{strong} rather than the \emph{weak} definition of forcing (see  \cite[Ch.13 1.8]{troelstravandalen:88ii}); that is,  a set of nodes is a cover/bar of a node p if it contains an entire \emph{level} over p. It is  not sufficient that every path through p passes through it. By a ``level $n$'' over p we mean all successors over $p$ at  hight $n$.

For enumerable first-order \theory\ we specialize \ref{proposition: modified joyals theorem for coherent with C coverage}/\ref{corollary: generalized beth model} to the effect that for  every $\theory^m$-model diagram $(D,F)$  in suitable \cat{M} there is a Beth model \alg{B} with root domain $D$ such that   $\alg{B}\Vdash \phi[\alg{d}/\alg{x}]\Leftrightarrow (D,F)\vDash P_{\phi}[\alg{d}/\alg{x}]$ for all first-order $\phi$. This yields a Beth$^{\star}$ completeness theorem for \theory.

 Specifically,  let \theory\ be a enumerable\ first-order theory over a signature $\Sigma$ and $\theory^m$ its regular Morleyization. In this section, we assume that  the sequent $\top\vdash\fins{x} x=x$ is an axiom of \theory. We refer to theories having this sequent as an axiom as \emph{habitative} Let \cat{M} be the subcategory of $\operatorname{MDiag}_b(\theory^m)$ consisting of diagrams $(D,F)$ of the following form. The domain $D$ is a semi-decidable subset of \thry{N}, coding a bounded relation from \thry{N} to \thry{N}. $D$ comes equipped with the least upper bound.  Denote by $f_D$ the function $\funktor{f_D}{\thry{N}}{2^{\thry{N}}}$ such that $D=\cterm{n\in \thry{N}}{\fins{m\in\thry{N}}f(n)(m)=1}$. The set of facts $F$ is, similarly, (coded as) a semi-decidable subset of \thry{N} (with function $f_F$). Then, straightforwardly, \cat{M} is chase-complete (for $\theory^m$) with a weakly reflective chase functor. We can assume that if $(D,F)$ is a diagram and $(D',F')=\operatorname{Ch}(D,F)$  then $f_D\leq f_{D'}$ and $f_F\leq f_{F'}$. We fix an enumeration $g$ of (codes of) all possible facts of the form \pair{\syntob{\alg{x}}{P_{\phi\vee\psi}},\alg{d}},  where every such fact is revisited an infinite number of times.   

For any $(D,F)\in\cat{M}$ we define a binary tree $\cat{T}_{D,F}$ over   $(D,F)$, where $\cat{T}_{D,F}$ occurs as a subcategory of \cat{M}, as follows. The diagram $(D,F)$ is the root. At node $(D',F')$ at level $n$, if $g(n)$ is, say, $\pair{\syntob{\alg{x}}{P_{\phi\vee\psi}},\alg{d}}$ and $\fins{m\leq n}f_{F'}(g(n),m)=1$ then $(D',F')$ has the children $\operatorname{Ch}(D',F'\cup\{P_{\phi}[\alg{d}]\})$ and $\operatorname{Ch}(D',F'\cup\{P_{\psi}[\alg{d}]\})$. Else the children of  $(D',F')$  are both $(D',F')$  itself. In the former case, we say that $\pair{\syntob{\alg{x}}{P_{\phi\vee\psi}},\alg{d}}$ \emph{is chased}.

$\cat{T}_{D,F}$ becomes a fallible Beth structure for $\Sigma$ by setting $P=\cat{T}_{D,F}$, and letting the functor $\funktor{\mathfrak{D}}{\cat{T}_{D,F}}{\topo{F}}$ send a $\Sigma^m$-diagram $S$ to its $\Sigma$ reduct extended with $\sem{\bot}^S:=\sem{P_{\bot}}^S$.
%
%  where the extension of an atomic (or Horn) formula $\phi$ in a node $(D',F')$ is set to be the extension of the formula $P_{\phi}$ in the $\theory^m$-model diagram $(D',F')$. Similarly, we define $(D',F')\Vdash \bot \Leftrightarrow (D',F')\vDash P_{\bot}$ 
%
It is clear that for any node $q$ of $\cat{T}_{D,F}$, any level of $\cat{T}_{D,F}$ above it  can be seen as a $\overline{C}$-cover of $q$ (and therefore also a $\overline{B}$-cover, since $\overline{B}=\overline{C}$ in this setting). We refer to it therefore also as ``a cover of $q$''. By the proof of Theorem \ref{proposition: modified joyals theorem for coherent with C coverage},  we have, therefore, that for any atomic formula $\FIC{\alg{x}}{P_{\phi}}$ of $\Sigma^m$, any node $q$ in $\cat{T}_{D,F}$, and any elements $\alg{d}$ of $q$, if $P_{\phi}[\alg{d}/\alg{x}]$ is true on a cover of $q$ then $q\vDash P_{\phi}[\alg{d}/\alg{x}]$.     

\begin{lemma}\label{lemma: nodes force the interpretation} Let $(D,F)\in \cat{M}$ and consider the Beth-structure $\cat{T}_{D,F}$. Let $p$ be a node in $\cat{T}_{D,F}$,
let $\FIC{\alg{x}}{\phi}$ be a first-order formula over $\Sigma$, with $P_{\phi}$ the corresponding  atomic formula in $\Sigma^m$, and let $\alg{d}\in p$. Then 
\[p\Vdash \phi[\alg{d}/\alg{x}]\Leftrightarrow p\vDash P_{\phi}[\alg{d}/\alg{x}]\]  
\begin{proof}
By induction on $\phi$, as follows.

$\phi$ atomic or $\phi = \bot$: by the remark immediately preceding this lemma.

$\phi\equiv \varphi\wedge \psi$ or $\phi=\top$: immediate.

$\phi= \varphi\vee \psi$: Assume $p\Vdash (\varphi\vee \psi)[\alg{d}]$. Then there exists a cover of $p$ such that for all $q$ in the cover either $q\Vdash \varphi[\alg{d}]$ or $q\Vdash \psi[\alg{d}]$. By induction hypothesis, then, either $q\vDash P_{\varphi}[\alg{d}]$ or $q\vDash P_{\psi}[\alg{d}]$. As $\theory^m$ proves e.g.\ $P_{\varphi}\vdash_{\alg{x}}P_{\varphi\vee\psi}$ therefore $q\vDash P_{\varphi\vee\psi}[\alg{d}]$. Whence $p\vDash P_{\varphi\vee\psi}[\alg{d}]$. 

\noindent Conversely, assume that $p\vDash P_{\varphi\vee\psi}[\alg{d}]$. Then for all $q\geq p$, we have $q\vDash P_{\varphi\vee\psi}[\alg{d}]$. Therefore, there exists a cover of $p$, say at level $n$, such that on that cover $P_{\varphi\vee\psi}[\alg{d}]$ is chased. Then for all $q\geq p$ at level $n+1$ we have that either $q\vDash P_{\varphi}[\alg{d}]$ or $q\vDash P_{\psi}[\alg{d}]$. Whence, by induction hypothesis, either $q\Vdash {\varphi}[\alg{d}]$ or $q\Vdash {\psi}[\alg{d}]$, and so $p\Vdash (\varphi\vee \psi)[\alg{d}]$. 

$\phi=\fins{y}\varphi$: We have that $p\Vdash \fins{y}\varphi[\alg{d}]$ iff $q\Vdash \varphi[a_q,\alg{d}]$ on a cover iff $q\vDash P_{\varphi}[a_q,\alg{d}]$ on a cover iff $q\vDash \fins{y}P_{\varphi}[\alg{d}]$ on a cover iff $q\vDash P_{\fins{x}\varphi}[\alg{d}]$ on a cover iff  $p\vDash P_{\fins{y}\varphi}[\alg{d}]$ .
 
$\phi= \varphi\rightarrow \psi$: Assume that $p\vDash P_{\varphi\rightarrow \psi}[\alg{d}]$. Let $q\geq p$ and assume that $q\Vdash \varphi[\alg{d}]$. By induction hypothesis  $q\vDash P_{\varphi}[\alg{d}]$. Now, $\theory^m$ proves the sequent $P_{\varphi}\wedge P_{\varphi\rightarrow \psi}\vdash_{\alg{x}}P_{\psi}$. Whence $q\Vdash P_{\psi}[\alg{d}]$. And so applying the induction hypothesis again, $q\Vdash {\psi}[\alg{d}]$. Hence $p\Vdash (\varphi\rightarrow \psi)[\alg{d}]$. 

\noindent For the converse direction, observe first that for $(D,F)\in \cat{M}$ we have that $(D,F)\vDash P_{\varphi\rightarrow \psi}[\alg{d}]$ iff  $\operatorname{Ch}(D,F\cup \{P_{\varphi}[\alg{d}]\})\vDash P_{\psi}[\alg{d}]$: for the right-to-left direction, the right hand side  implies that there is a $\theory^m$-provable sequent $P_{\chi}\wedge P_{\phi}\vdash_{\alg{x}} P_{\psi}$  such that $(D,F)\vDash P_{\chi}[\alg{d}/\alg{x}]$. Whence \theory\ proves ${\chi}\wedge {\phi}\vdash_{\alg{x}} {\psi}$ and therefore $\chi\vdash_{\alg{x}} \phi\rightarrow \psi$, so that $\theory^m$ proves  $P_{\chi}\vdash_{\alg{x}} P_{\phi\rightarrow \psi}$.  
Now, assume that $p\Vdash (\varphi\rightarrow \psi)[\alg{d}]$. We have $p\vDash P_{\top\vee \varphi}[\alg{d}]$, so there exists a level $n$ where it is chased. For a node $q=(D',F')$ on level $n$, the right child is therefore $q'=\operatorname{Ch}(D,F\cup \{P_{\varphi}(\alg{d})\})$. By induction hypothesis, $q'\vDash P_{\varphi}[\alg{d}]$ implies that $q'\Vdash {\varphi}[\alg{d}]$. Thus since $q'\geq p$ we have by assumption that $q'\Vdash {\psi}[\alg{d}]$, so $q'\vDash P_{\psi}[\alg{d}]$. So, by the observation, $q\vDash P_{\varphi\rightarrow\psi}[\alg{d}]$.  With $P_{\varphi\rightarrow\psi}[\alg{d}]$ true on a cover of $p$ we have, then, that $p\vDash P_{\varphi\rightarrow\psi}[\alg{d}]$.

 $\phi=\alle{y}\varphi$: Assume that $p\vDash P_{\alle{y}\varphi}[\alg{d}]$. Then for all $q\geq p$ and $a\in q$ we have that $q\vDash P_{\varphi}[a,\alg{d}]$, by applying the axiom $P_{\alle{y}\varphi}\vdash_{y,\alg{x}}P_{\varphi}$ of $\theory^m$. So, by induction hypothesis, $q\Vdash {\varphi}[a,\alg{d}]$. Thus $p\Vdash {\alle{x}\varphi}[\alg{d}]$.

\noindent For the converse direction, observe first, similar to the case of the conditional above, that that for $(D,F)\in \cat{M}$ we have that $(D,F)\vDash P_{\alle{x}\varphi}[\alg{d}]$ iff $\operatorname{Ch}(D+1,F)\vDash P_{\varphi}[e,\alg{d}]$, for all elements $e$ of $\operatorname{Ch}(D+1,F)$ (specifically the new element of $D+1$). Since \theory\ is habitative, this implies that $(D,F)\vDash P_{\alle{x}\varphi}[\alg{d}]$ iff $\operatorname{Ch}(D,F)\vDash P_{\varphi}[e,\alg{d}]$, for all elements $e$ of $\operatorname{Ch}(D,F)$ (specifically any fresh element $e$ added by an application of the axiom $\top\vdash\fins{x}x=x$).  Assume, then, that    $p\Vdash {\alle{x}\varphi}[\alg{d}]$. We have $p\vDash P_{\top\vee\top}$, so there exist a level $n$ on which $P_{\top\vee\top}$ is chased. For $(D',F')=q\geq p$ on level $n$, either child is $\operatorname{Ch}(D',F')$ and $\operatorname{Ch}(D',F')\Vdash \varphi[e,\alg{d}]$, for all elements $e$ in $(D',F')$. By induction hypothesis, $\operatorname{Ch}(D',F')\vDash P_{\varphi}[e,\alg{d}]$, so $q=(D',F')\vDash P_{\alle{x}\varphi}[\alg{d}]$. With   $ P_{\alle{x}\varphi}[\alg{d}]$ thus true on a cover, we conclude $p\vDash P_{\alle{x}\varphi}[\alg{d}]$.
\end{proof} 
\end{lemma}
%
%
%Say that a Beth structure for a first-order language \cat{L} is a \emph{Beth$^{\star}$ structure} if covers are only relevant for the forcing of  disjunctions; that is to say, for any node $p$, and elements $\alg{d}\in p$, the following holds (with the definition of the forcing relation being the usual, strong, one)
%\begin{enumerate}[(i)]
%\item for all atomic formulas $\phi\in \cat{L}$, it is the case that $p\Vdash\phi(\alg{d})$ iff $p\vDash\phi(\alg{d})$, and
%\item for all formulas of the form $\fins{x}\phi\in \cat{L}$, it is the case that  $p\Vdash\fins{x}\phi(\alg{d})$ iff there exists $a\in p$ such that $p\Vdash\phi(a,\alg{d})$.  
%\end{enumerate} 
%
It is clear that given a Beth$^{\star}$ model of a first order  \theory, the nodes, and the root in particular, model $\theory^m$. We can now state the converse 
\begin{theorem}
Let \theory\ be a habitative enumerable first-order theory. Let $\theory^m$ be its regular Morleyization. For every enumerable Tarski model \alg{M} of $\theory^m$ there exists a fallible Beth$^{\star}$ model \alg{B} of \theory\ such that the domain of the root $r$ is the domain of $\alg{M}$, and such that for all $\alg{m}\in \alg{M}^{l(\alg{m})}$ and \FIC{\alg{x}}{\phi} first-order, 
\[r\Vdash \phi[\alg{m}/\alg{x}]\Leftrightarrow \alg{M}\vDash P_{\phi}[\alg{m}/\alg{x}]\]  
%\begin{proof}
%\end{proof}
\end{theorem}
The following can then be seen as a constructive version of the completeness theorem of \cite{gabbay:77}. 
\begin{corollary}\label{theorem: Beth completeness for s.d.} Let \theory\ be a habitative enumerable\ first-order theory. Then \theory\ is complete with respect to fallible Beth$^{\star}$ models.

%\begin{proof}
%By Theorem \ref{proposition: modified joyals theorem for coherent with C coverage} and Lemma \ref{lemma: nodes force the %interpretation}.
%\end{proof} 
\end{corollary}

\section*{Acknowledgements}

We wish to particularly thank  H\aa{}kon R. Gylterud,  Hugo Herbelin, Panagis Karazeris,  Peter L.\ Lumsdaine, and, especially, Erik Palmgren, for discussions, pertinent observations, and feedback.

%\nocite{jech:08}
%\nocite{myhill:73}

\bibliographystyle{plain}
\bibliography{bibliografiEv}

\newpage

\end{document}